\documentclass[a4paper, 11pt, final]{amsart}

\makeatletter

\numberwithin{equation}{section}
\numberwithin{figure}{section}
\usepackage[a4paper]{geometry}
\usepackage{amsmath,amsfonts,amssymb,amsthm,epsfig}
\usepackage{amsmath}
\usepackage{color,xcolor}
\usepackage[final,linkcolor = blue,citecolor = blue,colorlinks=true]{hyperref}
\usepackage[leqno]{amsmath}
\usepackage[numbers,sort&compress]{natbib}
\usepackage{amscd}
\usepackage{mathrsfs}
\usepackage{}
\usepackage[all]{xy}
\usepackage[OT1]{fontenc}
\usepackage[latin1]{inputenc}
\usepackage{amssymb,latexsym}
\usepackage{type1cm,ae}
\usepackage[notref,notcite]{showkeys}
\usepackage{graphicx}
\usepackage{xspace}
\usepackage{setspace}
\usepackage{enumerate}
\usepackage[english]{babel}
\usepackage{bbm}
\usepackage{esint}
\usepackage{soul}

\usepackage{pifont}

\usepackage{graphicx}

\soulregister{\cite}7
\soulregister{\citep}7
\soulregister{\citet}7
\soulregister{\ref}7
\soulregister{\pageref}7
\sethlcolor{yellow}

\theoremstyle{definition}
\newtheorem{theorem}{Theorem}[section]
\newtheorem{proposition}[theorem]{Proposition}

\newtheorem{lemma}[theorem]{Lemma}

\newtheorem{definition}[theorem]{Definition}
\newtheorem{remark}[theorem]{Remark}

\numberwithin{equation}{section}

\newcommand*{\Wert}{\mathord{\mbox{|\kern-1.5pt|\kern-1.5pt|}}}

\newcommand{\N}{\mathbb{N}}	
\newcommand{\C}{\mathbb{C}}	
\newcommand{\R}{\mathbb{R}}	
\newcommand{\T}{\mathbb{T}}
\newcommand{\E}{\mathbb{E}}	
\newcommand{\A}{\mathbb{A}}	
\renewcommand{\P}{\mathbb{P}}	

\newcommand{\de}{\partial}		

\DeclareMathOperator{\supp}{supp}

\def\rr{{\Bbb R}}
\def\rz{{{\rr}^n}}

\def\cc{{\Bbb C}}
\def\de{{\Bbb D}}

\def\E{\mathbf{E}}

\def\L{\mathcal{L}}
\def\C{\mathcal{C}}
\def\M{\mathcal{M}}
\def\Ez{\mathcal{E}}
\def\Fz{\mathcal{F}}
\def\A{\mathcal{A}}

\def\W{\mathcal{W}}
\def\D{\mathcal{D}}

\def\P{\mathcal{P}}
\def\Q{\mathcal{Q}}
\def\H{\mathcal{H}}
\def\T{\mathcal{T}}
\def\R{\mathcal{R}}

\def\rdm{\rr^{n+1}}
\def\rdt{\rr_+^{n+2}}
\def\fz{\infty}

\def\supp{{\rm{\ supp\ }}}

\def\ez{\epsilon}

\def\supp{{\rm supp}}

\def\l{\left}
\def\r{\right}

\def\XXint#1#2#3{{\setbox0=\hbox{$#1{#2#3}{\int}$}
		\vcenter{\hbox{$#2#3$}}\kern-.5\wd0}}

\makeatother

\usepackage{babel}

\begin{document}
\raggedbottom
\allowdisplaybreaks

\title[The Kato problem for weighted elliptic and parabolic operators of higher order]{The Kato problem for weighted elliptic and parabolic operators of higher order}
	
	\author[Guoming Zhang]{Guoming Zhang}
	
	\address[Guoming Zhang]{College of Mathematics and System Science, Shandong University of Science and Technology, Qingdao, 266590, Shandong, People's Republic of China}
	\email{ zhangguoming256@163.com}

	

	\thanks{$^*$Corresponding author: Guoming Zhang}
	\thanks{The author is supported by the National Natural Science Foundation of Shandong Province (No. ~ZR2023QA124).}
	
	\date{}
	
	\begin{abstract} 
	
	We solve the Kato square root problem for parabolic operators of arbitrary order $2m$ whose coefficients are allowed to depend on both space and time in a merely measurable way and possess boundedness and ellipticity controlled by a Muckenhoupt $A_{2}-$weight. Notably, the proof applies to the weighted Kato problem within an elliptic framework. 
	\end{abstract}

	\maketitle
	\tableofcontents
\section{Introduction}
	
We study the $2m$ order parabolic operators in form of \begin{equation}\label{eq: c001}\H :=\partial_{t} +\L:=\partial_{t} +(-1)^{m}\sum_{|\alpha|=|\beta|=m}w^{-1}\partial^{\alpha}(a_{\alpha, \beta}\partial^{\beta}),\end{equation} where $m\in \N^{+}$ and $w=w(x)$ is independent of $t$ and belongs to the spatial Muckenhoupt class $A_{2}(\rz, dx).$ The coefficients $a_{\alpha, \beta}:=\{a_{\alpha, \beta}\}$ are complex-valued, measurable and dependent on both space and time, satisfying the degenerate ellipticity condition in the G$\mathring{a}$rding sense: \begin{equation}\label{eq: a2.28}\mbox{Re} \int_{\rdm}a_{\alpha, \beta}(t, x)\partial^{\alpha}f(t, x)  \overline{\partial^{\beta}f(t, x)}dxdt \geq c_{1} \int_{\rdm}|\nabla^{m}f(t, x)|^{2}w(x)dxdt, \quad \forall f\in \E_{\upsilon}, \end{equation} and the degenerate boundedness condition:\begin{equation}\label{eq: a2.29}|\sum_{|\alpha|=|\beta|=m} a_{\alpha, \beta}(t, x)\xi_{\alpha }\overline{\zeta_{\beta} }|\leq c_{2}w(x)|\xi| |\zeta|, \quad\forall \; \xi_{\alpha }, \;\zeta_{\beta}\in \cc,\end{equation} for some positive constants $c_{1}, c_{2}$ and all $(t, x)\in \rr^{n+1}.$ Here, for any $\xi, \zeta\in (\cc)^{p},$ $\xi \cdot \overline{\zeta}:=\sum_{|\beta|=m}\xi_{\beta} \overline{\zeta_{\beta}}$ denotes the inner product on $(\cc)^{p}$ (bear in mind that $(\cc)^{p}=\cc^{n}$ when $m=1$). $\E_{\upsilon}:=\E_{\upsilon}(\rdm; \cc)$ is the weighted parabolic energy space of $m$ order, which contains all functions $u$ such that $u, \nabla^{m}u:=(\partial^{\alpha}u)_{|\alpha|=m}$ and $D_{t}^{1/2} u$ are in the weighted Lebesgue space $L^{2}_{\upsilon}:=L^{2}(\rdm, d\upsilon) (d\upsilon:=dwdt=w(x)dxdt).$ The rigorous definitions of these objects can be found in Section 2.

This work is devoted to the derivation of the Kato estimate for $\H,$ that is, \begin{equation}\label{eq: t1}\|\sqrt{\H} u\|_{L^{2}_{\upsilon}}\approx \|\nabla^{m}u\|_{L^{2}_{\upsilon}}+ \|D_{t}^{1/2}u\|_{L^{2}_{\upsilon}}\quad (u\in \E_{\upsilon}),\end{equation} where the implicit constant depends only on $n, m,$ the ellipticity constant and the $A_{2}-$weight constant of $w.$ More precisely, we aim to show:  

\begin{theorem}\label{theorem: a0}\; Any operator $\H$ given in \eqref{eq: c001}-\eqref{eq: a2.29} can be defined as a maximal accretive operator in $L^{2}_{\upsilon}$ via an accretive sesquilinear form with domain $\E_{\upsilon}.$ Moreover, the domain of its unique maximal accretive square root $\sqrt{\H}$ coincides with $\E_{\upsilon},$ and the Kato estimate \eqref{eq: t1} is satisfied.

\end{theorem}
In particular, as a byproduct of the proof of \eqref{eq: t1}, an elliptic version of Theorem \ref{theorem: a0} is also obtained. Let $\L$ be as in \eqref{eq: c001}, defined on $\rz,$ and $W^{m, 2}_{w}:=W^{m, 2}_{w}(\rz; \cc)$ be the usual weighted Sobolev space of order $m$.   

\begin{theorem}\label{theorem: a00}\; Given $w\in A_{2}(\rz, dx).$ Suppose that the coefficients of $\L$ satisfy \eqref{eq: a2.28} -\eqref{eq: a2.29} on $W^{m, 2}_{w}$ and $\rz.$ Then $\L$ can be defined as a maximal accretive operator in $L^{2}_{w}:=L^{2}(\rz, dw)$ via an accretive sesquilinear form with domain $W^{m, 2}_{w},$ moreover, the domain of its unique maximal accretive square root $\sqrt{\L}$ is also $W^{m, 2}_{w}$ and 

\begin{equation}\label{eq: t2}\|\sqrt{\L} u\|_{L^{2}_{w}}\approx \|\nabla^{m}u\|_{L^{2}_{w}} \quad (u\in W^{m, 2}_{w})\end{equation} holds with the implicit constant depending only on $n, m,$ the ellipticity constant and the $A_{2}$-weight constant of $w.$
\end{theorem}

The Kato problem for unweighted elliptic operators (i.e. $w\equiv 1$) first posed by Kato \cite{TK} in 1961 was finally solved in the remarkable paper \cite{AHLT} by Auscher, $et\; al.$ The techniques introduced in \cite{AHLT} are highly effective in extending and applying to other problems, especially the $L^{p}$ (higher-order) Kato problems \cite{AHMT1, AT, AUS} and the boundary value problems for elliptic and parabolic equations and systems \cite{AAA, AA0, AA, AEK, AM1, AA1, CN, EH, HCMP, HLMPJ, HLM, NK0}. In $A_{2}-$weighted case, the elliptic Kato problem was first resolved in \cite{CR1} by extending the main techniques \cite{AHLT} to the weighted setting, then re-discovered in \cite{AAR} by using a different approach (i.e. the Dirac operator framework). An interesting extension of the Kato problem for degenerate elliptic operators appeared in \cite{DMR} and the $L^{p}-$version of the weighted Kato problem was settled in \cite{YZ}. 


In another direction, Nystr\"{o}m \cite{NK} extended the techniques in \cite{AHLT} to the unweighted parabolic setting and proved the equivalence between the square function estimates and Theorem \ref{theorem: a0} when $m=1$ and $w\equiv 1$ under the additional assumption that the coefficients are $t-$independent. Subsequently, Auscher, Egert and Nystr\"{o}m \cite{AEN1} established the unweighted parabolic Kato estimate in the absence of the $t$-independence of the coefficients through the Dirac operator framework. In a recent work, Ataei, Egert and Nystr\"{o}m \cite{AMK} utilized the techniques developed in \cite{NK} to provide a significantly simplified proof for Theorem \ref{theorem: a0} in the case $m=1,$ assuming that the coefficients depend measurably on all variables. Inspired by \cite{AMK, CR1}, we generalize the techniques rooted in \cite{AHLT, NK} to the higher-order and weighted setting in order to build the Kato estimates \eqref{eq: t1}-\eqref{eq: t2} in this paper. 

Intuitively, Theorem \ref{theorem: a0} extends \cite[Theorem 1.1]{AMK} to higher order parabolic operators. We therefore follow the same strategy as in \cite{AMK} to prove Theorem \ref{theorem: a0} and maintain the novelty in the same paper, that is: combining measurable dependence of the coefficients on all variables with $A_{2}-$weighted degenearacy in space and avoiding the Dirac operator framework. Our proofs are unavoidably technical, relying on the results and techniques proved previously for second order operators. The main technical lemmas cited in our proof are organized in the Appendix (Section 7).  In Section 2 we introduce some basic notations and definitions and gather some essential properties about $A_{2}-$weights and the relevant weighted energy spaces. 
In Section 3, the maximal accretivity of the part of $\H$ in $L^{2}_{\upsilon}$ is proved. Building upon this result and a duality argument, the Kato inequality \eqref{eq: t1} is reduced to a quadratic estimate:  
\begin{equation}\label{eq: d1}\int_{0}^{\infty}\|\lambda^{m}\H (1+\lambda^{2m}\H)^{-1} f\|^{2}_{L^{2}_{\upsilon}}\frac{d \lambda}{\lambda}\lesssim \|\nabla^{m} f\|^{2}_{L^{2}_{\upsilon}}+\|D_{t}^{1/2} f\|^{2}_{L^{2}_{\upsilon}}\quad (\forall \;f\in \E_{\upsilon})\end{equation} by employing the bounded $H^{\infty}-$calculus for maximal accretive operators in \cite{H, M}. The remaining sections are devoted to the proof of \eqref{eq: d1} and the proof is divided into three parts. First, we develop bounds and off-diagonal estimates for the resolvent operators in Section 4. It is important to note that, in contrast to the second order case, the lower order derivatives $\partial^{\alpha} (|\alpha|<m)$ play a role in the proof of off-diagonal estimates for the resolvent operators, which are ultimately addressed by the weighted Sobolev interpolation inequality. 
The second step involves using a higher-order version of the weighted Littlewood-Paley theory, the adapted parabolic setup and the boundedness of the ``principal part approximation operator'' to attain the estimate: \begin{equation}\label{eq: s2}
|||\lambda^{m}\H \Ez_{\lambda} f-(-1)^{m}\sum_{|\alpha|=|\beta|=m}\lambda^{m}\Ez_{\lambda}w^{-1}\partial^{\alpha}(a_{\alpha, \beta}\A_{\lambda}\partial^{\beta}f)|||_{2, \upsilon}\lesssim  \|\nabla^{m} f\|^{2}_{L^{2}_{\upsilon}}+\|D_{t}^{1/2} f\|^{2}_{L^{2}_{\upsilon}},\end{equation} see the details in Section 5. In the third step, we first reduce the problem to showing that 
\begin{equation}\label{eq: d2}|(-1)^{m}(\sum_{|\alpha|=m}\lambda^{m}(1+\lambda^{2m}\H)^{-1}w^{-1}\partial^{\alpha}(a_{\alpha, \beta}))_{|\beta|=m}|^{2}\frac{d\upsilon d\lambda}{\lambda}\;\;\mbox{is a Carleson measure.}\end{equation} Then, we generalize the weighted $Tb-$procedure in \cite[Section 8]{AMK} to higher order parabolic operators and use it to prove \eqref{eq: d2}. We stress that the smart trick in \cite{AMK}, namely, separating time and space variables, is sufficiently exploited in pursuit of constructing a higher-order version of the weighted Littlewood-Paley theory and the weighted $Tb-$procedure. This trick also explains why the coefficients are allowed to depend on all variables.  

Extending Theorem \ref{theorem: a0} and Theorem \ref{theorem: a00} to systems works without difficulty. We also note that Theorem \ref{theorem: a00} generalizes some partial results of \cite{AHMT1} to the weighted setting. Based on this paper and \cite[Theorem 1.11]{AHMT1}, we can expect the weighted Kato estimates associated to non-homogeneous higher-order elliptic operators (systems) $$\sum_{|\alpha|, |\beta|\leq m}(-1)^{|\alpha|}w^{-1}\partial^{\alpha}(a_{\alpha, \beta}\partial^{\beta})$$ and their parabolic counterparts. Additionally, Theorem \ref{theorem: a00} serves as the starting point for establishing a higher-order version of \cite{DMR, YZ}. On the other hand, our results offer significant potential for boundary value problems generated by elliptic and parabolic operators (systems) of higher order in the weighted context. 

The organization of the paper is outlined above.

\section{Preliminaries}

We now rigorously define the notations introduced earlier and introduce additional symbols to state our results.

\noindent{\textbf{2.1. Notation.}} For any fixed $(t, x) \in \rr\times \rz,$ we call $||(t, x)||_{m}$ the parabolic norm of $(t, x)$ if it is the unique solution $\rho$ of the equation $$\frac{t^2}{\rho^{4m}}+\sum_{i=1}^n\frac{x_{i}^2}{\rho^2}=1.$$ If $m=1,$ the parabolic norm $||(t, x)||_{m}$ becomes the usual one as in \cite{AMK, AEK, AEN1, CN, NK, NK0}. Given a half-open cube $Q=Q_{r}(x):=(x-r/2, x+r/2]^{n}$ with sidelength $r$ and center $x$ in $\rz,$ and an interval $I=I_{r}(t):=(t-r^{2m}/2^{2m}, t+r^{2m}/2^{2m}],$ we use $\Delta:=\Delta_{r}(t, x):=I\times Q\subset \rdm$ to denote a parabolic cube in $\rdm$ of size $l(\Delta):=r.$ For any $\lambda>0,$ $\lambda \Delta=(\lambda^{2m}I)\times (\lambda Q):=(t-(\lambda r)^{2m}/2^{2m}, t+(\lambda r)^{2m}/2^{2m}]\times (x-(\lambda r)/2, x+(\lambda r)/2]^{n}$ is used to represent the dilation of $\Delta.$ In what follows, we use $1_{E}$ to denote the characteristic function of a set $E.$ 

\noindent{\textbf{2.2. Weights.}} A real-valued and non-negative function $w(x)$ defined on $\rz$ is said to belong to the Muckenhoupt class $A_{2}(\rz, dx)$ if \begin{equation}\label{eq: a2.61}[w]_{A_{2}}:=\sup_{Q} \l(\frac{1}{|Q|}\int_{Q}w(x)dx\r)\cdot \l(\frac{1}{|Q|}\int_{Q}w(x)^{-1}dx\r)<\infty,\end{equation} where the supremum is taken with respect to all cubes $Q\subset \rz.$ Using the weight function $w,$ we can define a measure on $\rz$ by simply setting $dw(x):=w(x)dx,$ which satisfies that there exist two constants $\eta\in (0, 1)$ and $\gamma>0,$ depending only on $n$ and $[w]_{A_{2}},$ such that, for any measurable subset $E\subset Q,$ 
 \begin{equation}\label{eq: a2.611}
 \gamma^{-1}\l(\frac{|E|}{|Q|}\r)^{\frac{1}{2\eta}}\leq \frac{w(E)}{w(Q)}\leq  \gamma\l(\frac{|E|}{|Q|}\r)^{2\eta}.\end{equation}  
Here, $|\cdot|$ is the Lebesgue measure on $\rz$ and $w(E):=\int_{E}dw.$ We also need the notation for weighted averages, expressed as follows: $$(f)_{E, w}:=\fint_{E}f(x)dw(x):=\fint_{w(E)}f(x)w(x)dx:=\frac{1}{w(E)}\int_{w(E)}f(x)w(x)dx,$$ where $0<w(E)<\infty$ and $f$ is a locally integrable on $\rz$ with respect to $dw.$ In particular, we abbreviate $(f)_{E}$ when $w\equiv 1.$ From \eqref{eq: a2.611}, it is easy to see that $dw$ is doubling, that is, there is a constant $D$ (called doubling constant), depending only on $n$ and $[w]_{A_{2}},$ such that 
 \begin{equation}\label{eq: h0}
 w(2Q)\leq D w(Q)\quad \mbox{for all cubes}\; Q\subset \rz.\end{equation} In addition, the measure $d w^{-1}$ is also doubling thanks to \eqref{eq: a2.61}, and the doubling constant is still denoted by $D$ for simplicity.
 Moreover, the weight function $w$ induces two measures on $\rdm,$ defined as 
 \begin{equation}\label{eq: h1}\begin{aligned}&\quad\quad d\upsilon:=d\upsilon(t, x):=w(x)dxdt\\
 &\quad\quad d\upsilon^{-1}:=d\upsilon^{-1}(t, x):=w^{-1}(x)dxdt.
\end{aligned}\end{equation} They are all doubling concerning parabolic cubes $\Delta\subset \rdm,$ with the same doubling constant $2^{2m}D$ (see also \cite[(2.4)]{AMK}).

\noindent{\textbf{2.3. Weighted energy spaces.}} As usual, we let $L^{2}_{w}:=L^{2}(\rz, dw)$ be the weighted Lebesgue space with its norm denoted by $\|\cdot\|_{L^{2}_{w}}:=\|\cdot\|_{2, w},$ then it is a Hilbert space. Below, we use $<\cdot, \cdot>_{2, w}$ to denote the inner product and $\|\cdot\|_{L^{2}_{w}\to L^{2}_{w}} (\|\cdot\|_{2\to 2, w})$ to denote the operator norm of linear operators on $L^{2}_{w}.$ Owing to \eqref{eq: a2.61}, we see \begin{equation}\label{eq: h5}L^{2}_{w}\subset L^{1}_{loc}(\rz).\end{equation} It is well-known that the class $C_{0}^{\infty}(\rz)$ of smooth and compactly supported functions serves as a dense subspace of $L^{2}_{w}.$ The same properties and notations apply to $L^{2}_{\upsilon}$ in $\rdm.$

\begin{definition}\label{defnition: q0}\;(Elliptic weighted Sobolev space of higher order).\; We let $W^{m, 2}_{w}:=W^{m, 2}(\rz, dw)$ \\
$:=W^{m, 2}_{w}(\rz, \cc)$ be the space of all $f\in L^{2}_{w}$ such that its distributional derivatives $\partial^{\alpha} f $ belong to $L^{2}_{w}$ for all $|\alpha|\leq m.$ The space $W^{m, 2}_{w}$ is equipped with the norm $$\|f\|_{W^{m, 2}_{w}}:=(\|f\|_{2, w}^{2}+\sum_{|\alpha|\leq m}\|\partial^{\alpha} f\|_{2, w}^{2})^{1/2}.$$
\end{definition} 
Clearly, $W^{m, 2}_{w}$ is a Hilbert space and $C_{0}^{\infty}(\rz)$ is dense in $W^{m, 2}_{w}$. In particular, using the weighted Sobolev interpolation inequality in \cite{GW}: \begin{equation}\label{eq: a2.32}\l(\int_{\rr^{n}}|\partial^{\gamma}v|^{2}w\r)\leq c_{0}(n, m, [w]_{A_{2}}) \l(\int_{\rr^{n}}|v|^{2}w\r)^{(1-\frac{|\gamma|}{m})}\l(\int_{\rr^{n}}|\nabla^{m}v|^{2}w\r)^{\frac{|\gamma|}{m}} \quad (\forall\; |\gamma|\leq m),\end{equation} we obtain, for any $f\in W^{m, 2}_{w},$
\begin{equation}\label{eq: h3}
\|f\|_{W^{m, 2}_{w}}\approx (\|f\|_{2, w}^{2}+\sum_{|\alpha|= m}\|\partial^{\alpha} f\|_{2, w}^{2})^{1/2}.\end{equation} 
For ease of presentation, we adopt the notations $\lesssim$ and $\approx.$ Specifically, for two positive constants $A, B$, the expression $A\lesssim B$ means that there exists a nonessential constant $C,$ depending only on $n, m, c_{1}, c_{2}$ and $[w]_{A_{2}},$ such that $A\leq C B.$ The notations $A \gtrsim B$ and $A\approx B$ should be interpreted similarly.

The operators $D_{t}^{1/2}$ and $H_{t}$ denote the half-order derivative and Hilbert transform in time through the Fourier symbols $|\tau|^{1/2}$ and $isgn(\tau),$ respectively. By \cite{NPV} (considering $f(\cdot, x)$ with $x\in \rz$ fixed), we have for $a.e. x\in \rz$ that 
\begin{equation}\label{eq: h2}\begin{aligned}
\|D_{t}^{1/2} f\|_{2, \upsilon}^{2}=c\int_{\rz}\int_{\rr}\int_{\rr}\frac{|f(t, x)-f(s, x)|^{2}}{|t-s|^{2}}dtdsdw(x)
\end{aligned}\end{equation} with the right-hand side being finite precisely when $D_{t}^{1/2} f\in L^{2}_{\upsilon}.$ 
We now introduce the higher-order parabolic weighted energy space in the parabolic setting. 

\begin{definition}\label{defnition: q1}\;(Parabolic weighted energy space of higher order).\; The space $\E_{\upsilon}:=\E_{\upsilon}(\rz, \cc)$\\$:=\E_{\upsilon}(\rdm, d\upsilon)$ consists of all functions $f\in L^{2}_{\upsilon}$ such that $D_{t}^{1/2} f$ and $\partial^{\alpha} f$ belong to $L^{2}_{w}$ for all $|\alpha|\leq m.$ The norm of $\E_{\upsilon}$ is defined by $$\|f\|_{\E_{\upsilon}}:=(\|f\|_{2, \upsilon}^{2}+\sum_{|\alpha|\leq m}\|\partial^{\alpha} f\|_{_{2, \upsilon}}^{2}+\| D_{t}^{1/2}f\|_{2, \upsilon}^{2})^{1/2}.$$
\end{definition} 
It is easy to show that $\E_{\upsilon}$ is a Hilbert space. Moreover, it follows from \eqref{eq: h3} that \begin{equation}\label{eq: h4}\begin{aligned}
\|f\|_{\E_{\upsilon}} \approx (\|f\|_{2, \upsilon}^{2}+\sum_{|\alpha|=m}\|\partial^{\alpha} f\|_{_{2, \upsilon}}^{2}+\| D_{t}^{1/2}f\|_{2, \upsilon}^{2})^{1/2}\end{aligned}\end{equation} holds for all $f\in \E_{\upsilon}.$ In the following,  we sometimes use the symbol $\de f:=(\nabla ^{m} f, D_{t}^{1/2}f)=(\partial^{\alpha}f,  D_{t}^{1/2}f)_{|\alpha|=m}$ when $f\in \E_{\upsilon}.$




\section{The higher order parabolic operator and the central estimate to \eqref{eq: t1}} 
In this section we first give a comprehensive discussion of how to interpret the parabolic operator \eqref{eq: c001} as either a bounded operator $\E_{\upsilon} \to \E_{\upsilon}^{*}$ or an unbounded operator on $L^{2}_{\upsilon}$ via a sesquilinear form. Hereafter, the superscript $^{*}$ stands for the dual.

Note that Lemma \ref{lemma: b1.5} implies that \begin{equation}\label{eq: b001}
\E_{\upsilon}\subset L^{2}_{\upsilon}\simeq (L^{2}_{\upsilon})^{*} \subset (\E_{\upsilon})^{*}.
\end{equation} By self-evident embeddings,  
 \begin{equation}\label{eq: b002}
D_{t}^{1/2}: \E_{\upsilon}\to L^{2}_{\upsilon},\quad \partial^{\alpha}: \E_{\upsilon}\to L^{2}_{\upsilon},\quad \forall \; |\alpha|\leq m,\end{equation} furthermore, in view of \eqref{eq: b001}, \begin{equation}\label{eq: b003}
D_{t}^{1/2}: L^{2}_{\upsilon}\to \E_{\upsilon}^{*},\quad w^{-1}\partial^{\alpha}(w\cdot ): L^{2}_{\upsilon}\to \E_{\upsilon}^{*},\quad \forall \; |\alpha|= m.\end{equation} If we split \begin{equation}\label{eq: c0}\partial_{t}=D_{t}^{1/2}H_{t}D_{t}^{1/2},\end{equation} it then follows from \eqref{eq: b002} and \eqref{eq: b003} that $$\partial_{t}: \E_{\upsilon}\to \E_{\upsilon}^{*}.$$ Consequently, $\H$ can be defined as a bounded operatior $\E_{\upsilon}\to \E_{\upsilon}^{*}$ via \begin{equation}\label{eq: b004}<\H u, \phi>:=\int_{\rdm} H_{t}D_{t}^{1/2}u \cdot\overline{D_{t}^{1/2}\phi}+\sum_{|\alpha|=|\beta|=m}w^{-1}a_{\alpha, \beta}\partial^{\beta}u \cdot\overline{\partial^{\alpha}\phi}d\upsilon \quad \forall\; u, \phi\in \E_{\upsilon}.\end{equation}

Due to \eqref{eq: b001}, we are justified in considering the maximal restriction of $\H$ to an operator in $L^{2}_{\upsilon},$ called the part of $\H$ in $L^{2}_{\upsilon},$ with domain $\D(\H):=\{u\in \E_{\upsilon}: \H u\in L^{2}_{\upsilon}\}.$ Obviously, $$<\H u, v>=\int_{\rdm} \H u\overline{\phi}$$ holds for all $\phi\in \E_{\upsilon}$ and $u\in \D(\H),$ which implies that the part of $\H$ in $L^{2}_{\upsilon}$ gives a meaning to the formal expression \eqref{eq: c001} by applying a formal integration by parts in \eqref{eq: b004}. In particular, it follows from \eqref{eq: b002} and \eqref{eq: b003} again that 

\begin{equation}\label{eq: b005}\H= D_{t}^{1/2}H_{t}D_{t}^{1/2}+\sum_{|\alpha|=|\beta|=m}(-1)^{|\alpha|}w^{-1}\partial^{\alpha}(w)(w^{-1}a_{\alpha, \beta}\partial^{\beta}).\end{equation} The equation \eqref{eq: b005} plays a key role in our proof.
 
\noindent{\textbf{3.1. Maximal accretivity.}} After introducing $\H$ as a bounded operator $\E_{\upsilon}\to \E_{\upsilon}^{*}$ in \eqref{eq: b004}, we need to further demonstrate that it is a maximal accretive operator on the Hilbert space $L^{2}_{\upsilon}$ in order to run the functional calculus for sectorial operators \cite{H, M}. Indeed, we have the following lemma.
 
\begin{lemma}\label{lemma: a12.6}\; Let $\sigma\in \cc$ with $\mbox{Re}\; \sigma>0.$  The following assertions are true.

$(i)\quad$ For each $f\in \E_{\upsilon}^{*},$  there exists a unique $u\in \E_{\upsilon}$ such that $(\sigma+\H)u=f.$ Moreover, $$\|u\|_{\E_{\upsilon}}\leq C \max\{\frac{c_{2}+1}{c_{1}}, \frac{|\mbox{Im} \;\sigma|+1}{\mbox{Re}\; \sigma}\}\|f\|_{\E_{\upsilon}^{*}},$$ and 
 for all $ \phi \in \E_{\upsilon},$ \begin{equation}\label{eq: a2.27}\int_{\rr^{n+1}}\sigma u\overline{\phi}+\sum_{|\alpha|=|\beta|=m}w^{-1}a_{\alpha, \beta}\partial^{\beta}u \cdot\overline{\partial^{\alpha}\phi}+H_{t}D_{t}^{1/2}u \cdot\overline{D_{t}^{1/2}\phi}d\upsilon=f(\phi).\end{equation} 

$(ii)\quad$ If $f\in L^{2}_{\upsilon},$ then \begin{equation}\label{eq: b008}\|(\sigma+\H)^{-1}f\|_{2, \upsilon}\leq \frac{ \|f\|_{2, \upsilon}}{\mbox{Re}\; \sigma}.\end{equation} In particular, the part of $\H$ in $L_{\upsilon}^{2}$ is maximal accretive and $\D(\H)$ is dense in $\E_{\upsilon}.$

$(iii)\quad$ The adjoint $\H^{*}$ of $\H$ can be identified formally as the backward-in-time operator $$-\partial_{t}+\sum_{|\alpha|=|\beta|=m}(-1)^{|\beta|}w^{-1}(\partial^{\beta}w)(w^{-1}\overline{a_{\alpha, \beta}}\partial^{\alpha}),$$ and all the above results hold with $\H$ replaced by $\H^{*}.$
\end{lemma} 
{\it Proof.}\quad We begin with $(i).$ The proof of $(i)$ relies on the hidden coercivity of the parabolic sesquilinear form \eqref{eq: b004}, a property originally discovered by Kaplan \cite{K} and re-discovered several times in \cite{AME, DZ, HL, NK}. The case $m=1$ is stated in \cite[Lemma 4.1]{AMK}.


We define the sesquilinear form $B_{\delta, \sigma}:\E_{\upsilon}\times \E_{\upsilon}\to \cc$ via 
\begin{equation}\label{eq: a00}\begin{aligned}
B_{\delta, \sigma}(u, \phi):=\int_{\rdm}\sigma u\cdot&\overline{(I+\delta H_{t})\phi}+w^{-1}a_{\alpha, \beta}\partial^{\beta}u\cdot \overline{\partial^{\alpha}(I+\delta H_{t})\phi}\\
&\quad\quad+H_{t}D_{t}^{1/2}u\cdot\overline{D_{t}^{1/2}(I+\delta H_{t})\phi}d\upsilon, 
\end{aligned}\end{equation} where $\delta\in (0, 1)$ to be chosen later. By Plancherel's theorem, the Hilbert transform $H_{t}$ is isometric on $\E_{\upsilon},$ then $B_{\delta, \sigma}$ is bounded. Observing that $H_{t}$ is skew-adjoint, we then see \begin{equation}\label{eq: b007}\mbox{Re}\;\int_{\rdm}H_{t}h \cdot \overline{h}=0 \quad \mbox{for all}\; h\in L^{2}_{\upsilon}.\end{equation} Exploiting \eqref{eq: a2.28}-\eqref{eq: a2.29} and \eqref{eq: b007}, it is now standard to deduce that, for any $\sigma \in \cc,$ \begin{equation}\label{eq: a01}
\mbox{Re}\; B_{\delta, \sigma}(u, u)\geq \delta\|D_{t}^{1/2}u\|_{2, \upsilon}^{2}+(c_{1}-c_{2}\delta)\|\nabla^{m}u\|_{2, \upsilon}^{2}+(\mbox{Re}\;\sigma-\delta |\mbox{Im}\;\sigma|)\|u\|_{2, \upsilon}^{2}.\end{equation} If we restrict $\mbox{Re}\; \sigma>0$ and choose $\delta:=\min\{\frac{c_{1}}{1+c_{2}}, \frac{\mbox{Re}\;\sigma}{1+|\mbox{Im}\;\sigma|}\}$ in the above inequality, then by \eqref{eq: h4}, \begin{equation}\label{eq: b006}
\mbox{Re}\; B_{\delta, \sigma}(u, u)\gtrsim \min\{\frac{c_{1}}{1+c_{2}}, \frac{\mbox{Re}\;\sigma}{1+|\mbox{Im}\;\sigma|}\}\|u\|_{\E_{\upsilon}}^{2}
\end{equation}
With \eqref{eq: b006} in hand, an application of the Lax-Milgram lemma and the fact that $I+\delta H_{t}$ is isometric on $\E_{\upsilon}$ yields \eqref{eq: a2.27}. We complete the proof of $(i).$

We now proceed with the proof of $(ii).$ By $(i)$, for any $\sigma\in \cc$ with $\mbox{Re}\; \sigma>0,$ we know that $\sigma +\H: \D(\H)\to L_{\upsilon}^{2}$ is bijective. Letting $f\in L_{\upsilon}^{2},$ we can define $u:=(\sigma +\H)^{-1}f\in\D(\H)\subset \E_{\upsilon}.$ Utilizing \eqref{eq: a2.28} and \eqref{eq: b007}, we get 

\begin{equation*}\begin{aligned} \mbox{Re}\; \sigma \|u\|_{2, \upsilon}^{2}&\leq\mbox{Re}\; \int_{\rdm}(\sigma u\cdot \overline{u}+\sum_{|\alpha|=|\beta|=m}a_{\alpha, \beta}w^{-1}\partial^{\beta}u \cdot\overline{\partial^{\alpha}u}+H_{t}D_{t}^{1/2}u \cdot\overline{D_{t}^{1/2}u})\;d\upsilon\\ 
 &=\mbox{Re}\; <(\sigma+\H) u, u>_{2, \upsilon}=\mbox{Re}\; \int_{\rdm} f\cdot \overline{u}d\upsilon \leq \|f\|_{2, \upsilon}\|u\|_{2, \upsilon}.\end{aligned}\end{equation*} Thus \begin{equation}\label{eq: b008}\|(\sigma+\H)^{-1}f\|_{2, \upsilon}\leq \frac{ \|f\|_{2, \upsilon}}{\mbox{Re}\; \sigma}.\end{equation} This implies that the resolvent set of the part of $\H$ in $L_{\upsilon}^{2}$ is non-empty, so the part of $\H$ in $L_{\upsilon}^{2}$ is closed. An application of \cite[Proposition 2.1.1]{H} shows that $\D(\H)$ is dense in $L_{\upsilon}^{2}.$ As a result, the part of $\H$ is maximal accretive. 
 
It remains to prove that $\D(\H)$ is indeed dense in $\E_{\upsilon}.$ To this end, we use the sesquilinear form $B_{\delta, 1}$ from \eqref{eq: a00}, with $\delta$ chosen to ensure \eqref{eq: b006}. Given a $\phi\in \E_{\upsilon}$ being orthogonal to $\D(\H).$ Thanks to \eqref{eq: b006}, we can apply the Lax-Milgram lemma to find a unique $\omega\in \E_{\upsilon}$ such that $$ <u, \phi>_{\E_{\upsilon}}=B_{\delta, 1}(u, \omega)\quad \mbox{for all}\; u\in  \E_{\upsilon}.$$ Confining further $u\in \D(\H),$ it follows from \eqref{eq: a2.27} that $$0=<(1+\H)u, (I+\delta H_{t})\omega>_{2, \upsilon}.$$ Since $1 +\H: \D(\H)\to L_{\upsilon}^{2}$ is bijective, the latter estimate implies $$0=<h, (I+\delta H_{t})\omega>_{2, \upsilon}\quad  \mbox{for all}\; h\in  L_{\upsilon}^{2}.$$ Therefore, $(I+\delta H_{t})\omega=0.$ This yields $\omega=0,$ then $\phi=0.$ The proof of $(ii)$ is complete.

Eventually, substituting the sesquilinear form $B_{\delta, \sigma}(u, \phi)$ with $B^{*}_{\delta, \sigma}(u, \phi):=\overline{B_{\delta, \sigma}(\phi, u)}$ and repeating the above arguments, we conclude with $(iii).$

\hfill$\Box$ 

\noindent{\textbf{3.2. Reducing \eqref{eq: t1} to the quadratic estimate \eqref{eq: d1}.}} 

Since $\H$ is maximal accretive according to Lemma \ref{lemma: a12.6}, then it has a unique accretive square root $\sqrt{\H}$ defined by the functional calculus for sectorial operators (see \cite{H, M}), and the same is true for the adjoint $\H^{*}$ with $(\sqrt{\H})^{*}=\sqrt{\H^{*}}.$ This allows us to apply \cite[Theorem 5.2.6]{H} and write \begin{equation}\label{eq: a2.55}
\sqrt{\H}f =c(m)\int_{0}^{\infty}\lambda^{3m}\H^{2}(1+\lambda^{2m}\H)^{-3}f\frac{d\lambda}{\lambda}
\quad (f\in \D(\sqrt{\H})).\end{equation} It is clear that the integral is understood as an improper Riemann integral in $L^{2}_{\upsilon}.$ With \eqref{eq: a2.55} in hand, for any $h\in L^{2}_{\upsilon},$ it follows that $$|\langle \sqrt{\H}f, h \rangle_{2, \upsilon}|\lesssim |||\lambda^{m}\H (1+\lambda^{2m}\H)^{-1} f|||_{2, \upsilon}\cdot  |||\lambda^{2m}\H^{*} (1+\lambda^{2m}\H^{*})^{-2} h|||_{2, \upsilon},$$ where $$||| \cdot |||_{2, \upsilon}:= \l(\iint_{\rdt}| \cdot|^2\frac{d\upsilon d\lambda}{\lambda}\r)^{1/2}.$$ Since $\H^{*}$ is also maximal accretive in $L^{2}_{\upsilon}$ thanks to the conclusion $(iii)$ in Lemma \ref{lemma: a12.6},  we can use \cite[Theorem 7.1.7; Theorem 7.3.1]{H} to deduce 
$$ |||\lambda^{2m}\H^{*} (1+\lambda^{2m}\H^{*})^{-2} h|||_{2, \upsilon} \lesssim \|h\|_{2, \upsilon}.$$ Consequently, by taking the supremum over all $h,$ we obtain 
$$ \|\sqrt{\H} f\|_{2, \upsilon}\lesssim |||\lambda^{m}\H \Ez_{\lambda} f|||_{2, \upsilon}.$$ 

Suppose that  \begin{equation}\label{eq: a2.54} |||\lambda^{m}\H (1+\lambda^{2m}\H)^{-1} f|||_{2, \upsilon}\lesssim \|\nabla^{m} f\|_{2, \upsilon} +\|D_{t}^{1/2}f\|_{2, \upsilon}\quad (f\in \E_{\upsilon})\end{equation} holds for now. 
Then \begin{equation}\label{eq: a2.56}
\|\sqrt{\H} f\|_{2,  \upsilon}\lesssim \|\nabla^{m} f\|_{2,  \upsilon} + \|D_{t}^{1/2}f\|_{2,  \upsilon}\quad (f\in \D(\sqrt{\H})\cap \E_{ \upsilon}\supset \D(\H)).\end{equation} By Lemma \ref{lemma: a12.6} again, $\D(\H)$ is dense in $ \E_{ \upsilon}.$ Bear in mind that $\sqrt{\H}$ is closed in $L^{2}_{ \upsilon}.$ This implies that \eqref{eq: a2.56} is indeed true for all $f\in \E_{ \upsilon}.$ In veiw of the transform $f(t, x)\to f(-t, x)$ and the definition of $\H^{*},$ we see that \eqref{eq: a2.56} also holds with $\sqrt{\H}$ replaced by $\sqrt{\H^{*}},$ that is, \begin{equation}\label{eq: a2.57}
\|\sqrt{\H^{*}} h\|_{2, \upsilon}\lesssim \|\nabla^{m} h\|_{2, \upsilon} + \|D_{t}^{1/2}h\|_{2, \upsilon}\quad (h\in \E_{\upsilon}).\end{equation} On the other hand, by using \eqref{eq: a01} with $\sigma=0$ and letting $\delta$ small, we deduce that, for all $f\in \D(\H),$ 
\begin{equation*}\begin{aligned}
\|\nabla^{m} f\|_{2, \upsilon}^{2} + \|D_{t}^{1/2}f\|_{2, \upsilon}^{2}\lesssim |\langle \H f, (I+\delta H_{t}) f \rangle_{2, \upsilon}|\lesssim \|\sqrt{\H} f\|_{2, \upsilon} \|\sqrt{\H^{*}}(I+\delta H_{t}) f\|_{2, \upsilon}.\end{aligned}
\end{equation*} This inequality together with \eqref{eq: a2.57} yield \begin{equation}\label{eq: a2.58}\|\nabla^{m} f\|_{2, \mu} + \|D_{t}^{1/2}f\|_{2, \upsilon} \lesssim \|\sqrt{\H} f\|_{2, \upsilon}\quad (f\in \D(\H)).\end{equation} Invoking \cite[Proposition 3.1.1(h)]{H} and the fact that $\D(\H)$ is dense in $\D(\sqrt{\H})$ with respect to the graph norm, it follows immediately that \eqref{eq: a2.58} holds for all $f\in \D(\sqrt{\H})$ and $\D(\sqrt{\H})=\E_{\upsilon}.$ Therefore, concatenating \eqref{eq: a2.56} and \eqref{eq: a2.58}, we successfully reduce Theorem \ref{theorem: a0} to the quadratic estimate \eqref{eq: a2.54}. For this reduction, one can also refer to the argument in \cite[Section 6]{AMK}.

Clearly, to obtain \eqref{eq: a2.54}, the norm estimates for the resolvent operator $(1+\lambda^{2m}\H)^{-1}$ are required. This motivates the next section. 

\section{Proof of \eqref{eq: a2.54}: Part I}

Set $\Ez_{\lambda}:=(1+\lambda^{2m}\H)^{-1}$ and $\Ez^{*}_{\lambda}:=(1+\lambda^{2m}\H^{*})^{-1}.$ Imitating the proof of \cite[Lemma 4.3]{AMK}, we have the following lemma.   
 
\begin{lemma}\label{lemma: a12.7}\; For all $\lambda>0$ and $f\in L^{2}_{\upsilon},$ \begin{equation*}\begin{aligned}&(i)\quad \|\Ez_{\lambda}f\|_{2, \upsilon}+\|\lambda^{m}\de \Ez_{\lambda}f\|_{2, \upsilon} \lesssim \|f\|_{2, \upsilon}\\ 
&(ii)\quad \|\lambda^{m}\Ez_{\lambda}D_{t}^{1/2} f\|_{2, \upsilon}+\|\lambda^{2 m}\de \Ez_{\lambda} D_{t}^{1/2}f\|_{2, \upsilon}\lesssim \|f\|_{2, \upsilon}\\ 
&(iii)\quad\|\lambda^{m}\Ez_{\lambda}w^{-1}(\nabla^{m} (wf))\|_{2, \upsilon}+\|\lambda^{2 m}\de \Ez_{\lambda} w^{-1}(\nabla^{m} (wf))\|_{2, \upsilon}\lesssim \|f\|_{2, \upsilon}.\end{aligned}\end{equation*}
The same estimates hold for $\Ez^{*}_{\lambda}.$ 
\end{lemma} 
{\it Proof.}\quad To achieve $(i),$ we set $u:=(\lambda^{-2m}+\H)^{-1}f.$ Then, $\lambda^{-2m}u=\Ez_{\lambda}f,$ and from \eqref{eq: b008}, it follows that $$\|\Ez_{\lambda}f\|_{2, \upsilon} \lesssim \|f\|_{2, \upsilon}.$$
Recalling the definition of sesquilinear form $B_{\delta, \lambda^{-2m}}$ in \eqref{eq: a00}, and applying the Lax-Milgram lemma, we conclude by choosing $\delta=\frac{c_{1}}{2 c_{2}}$ that $$\mbox{Re}\; B_{\delta,  \lambda^{-2m}}(u, u)\geq \delta(\|D_{t}^{1/2}u\|_{2, \upsilon}^{2}+\|\nabla^{m}u\|_{2, \upsilon}^{2})+ \lambda^{-2m}\|u\|_{2, \upsilon}^{2}$$ and
\begin{equation}\label{eq: b009}B_{\delta, \lambda^{-2m}}(u, u)=<f, (I+\delta H_{t})u>_{2, \upsilon}.\end{equation}  As a consequence, $$\|\de u\|_{2, \upsilon}\lesssim \lambda^{m} \|f \|_{2, \upsilon}. $$ This proves $(i).$ Since $\H^{*}$ shares the same structure with $\H$ from the point of view of sesquilinear form, the above arguments apply to $\Ez^{*}_{\lambda}.$ So $(i)$ holds for $\Ez^{*}_{\lambda}$ in place of $\Ez_{\lambda}.$

By duality in $L^{2}_{\upsilon}$ and the conclusion $(i)$ for $\Ez^{*}_{\lambda},$ we can derive $$\|\lambda^{m}\Ez_{\lambda}D_{t}^{1/2} f\|_{2, \upsilon}+\|\lambda^{m}\Ez_{\lambda}w^{-1}(\nabla^{m} (wf))\|_{2, \upsilon}\lesssim \|f\|_{2, \upsilon}.$$ 

We now consider $\lambda^{2 m}\de \Ez_{\lambda} D_{t}^{1/2}f.$ Note that $D_{t}^{1/2} f \in \E_{\upsilon}^{*}$ by \eqref{eq: b003}. Using Lemma \ref{lemma: a12.6}, we see $u:=(\lambda^{-2m}+\H)^{-1} D_{t}^{1/2}f\in \E_{\upsilon}.$
Plugging it into \eqref{eq: b009} we reach $$B_{\delta, \lambda^{-2m}}(u, u)=<f, D_{t}^{1/2}(I+\delta H_{t})u>_{2, \upsilon},$$ which yields $$\|\de u\|_{2, \upsilon}\lesssim  \|f \|_{2, \upsilon}. $$ Due to \eqref{eq: b003} again, we can replace 
$D_{t}^{1/2} f $ by $w^{-1}\nabla^{m} (wf)$ in the above derivation, which leads to $$\|\lambda^{2 m}\de \Ez_{\lambda}w^{-1}\nabla^{m} (wf)\|_{2, \upsilon}\lesssim \|f\|_{2, \upsilon}.$$ This ends the proof.

\hfill$\Box$ 

Given two measurable subsets $E, F$ in $\rdm,$ we define their parabolic distance by $$d(E, F):=\inf\{\|(t-s, x-y)\|_{m}: (t, x)\in F, (s, y)\in E\}.$$ As mentioned in \cite[Lemma 4.4]{AMK}, we cannot expect off-diagonal estimates for the non-local operator $D_{t}^{1/2}.$ Fortunately, the off-diagonal estimates estimates involving only the spatial derivatives are sufficient for our purposes, as demonstrated in \cite{AMK}. 



\begin{lemma}\label{lemma: a2.5}\;  Assume that $E, F$ are measurable subsets of $\rdm,$ and let $d:=d(E, F),$ then \begin{equation}\label{eq: a2.26}\int_{F}|\Ez_{\lambda}f |^{2}+|\lambda^{m} \nabla^{m}\Ez_{\lambda}f |^{2}d\upsilon\lesssim e^{-\frac{ d}{c\lambda}} \int_{E}|f |^{2}d\upsilon\end{equation} and \begin{equation}\label{eq: a2.261}\int_{F}|\lambda^{m} \Ez_{\lambda}[w^{-1}\partial^{\alpha}(wf)] |^{2}d\upsilon\lesssim e^{-\frac{ d}{c\lambda}} \int_{E}|f |^{2}d\upsilon, \quad \forall\; |\alpha|=m,\end{equation} hold for some constant $c,$ depending only on $n, m, c_{1}, c_{2}$ and $[w]_{A_{2}},$ and all $f\in L_{\upsilon}^{2}$ with $\supp f\subset E.$ The same statements are also true with $\Ez_{\lambda}$ replaced by $\Ez_{\lambda}^{*}.$
 \end{lemma} 
{\it Proof.}\quad
We follow the strategy of \cite[Lemma 4.4]{AMK} to prove \eqref{eq: a2.26}. First. we let $$ \Delta:=\frac{\kappa d}{\lambda}$$ with $0<\kappa<1$ to be chosen later. Clearly, we can assume $\Delta\geq 1,$ otherwise \eqref{eq: a2.26} is trivial. Second, we pick $\tilde{\eta}\in C^{\infty}(\rr^{n+1})$ such that $\tilde{\eta}=1$ on $F$ and $\tilde{\eta}=0$ on $E,$ also satisfying, for all $|\alpha|\leq m,$
\begin{equation}\label{eq: a2.30}|\partial^{\alpha}\tilde{\eta}|\lesssim d^{-|\alpha|}, \quad |\partial_{t}\tilde{\eta}|\lesssim d^{-2m}.\end{equation} Setting $\eta=e^{\Delta\tilde{\eta}}-1,$ we then have $$u:= \Ez_{\lambda}f \in \E_{\upsilon}\quad \mbox{and}\quad v:=u\eta^{2}\in \E_{\upsilon}\mbox{ thanks to Lemma \ref{lemma: b1.5} }.$$ 

By the density of $C_{0}^{\infty}(\rr^{n+1})$ in $\E_{\upsilon}$ (Lemma \ref{lemma: b1.5}) and the argument at the top of \cite[Page 11]{AMK}, we conclude that $$\mbox{Re}\; \int_{\rr^{n+1}}H_{t}D_{t}^{1/2}u \cdot\overline{D_{t}^{1/2}v}d\upsilon=-\frac{1}{2}\int_{\rr^{n+1}}|u|^{2}\partial_{t}(\eta^{2})d\upsilon.$$ By this, \eqref{eq: a2.27} with $\sigma =\lambda^{-2m}$ and $\supp f\subset E$, it follows that \begin{equation}\label{eq: a2.31}\begin{aligned}\int|u|^{2}\eta^{2}+\lambda^{2m} \mbox{Re}\;\sum_{|\alpha|=|\beta|=m}w^{-1}a_{\alpha, \beta}\partial^{\beta}u \cdot\overline{\partial^{\alpha}(u\eta^{2})}d\upsilon=\frac{\lambda^{2m}}{2}\int_{\rr^{n+1}}|u|^{2}\partial_{t}(\eta^{2})d\upsilon.\end{aligned}\end{equation} In virtue of the definition of $\eta,$ we can rewrite \eqref{eq: a2.31} as \begin{equation}\label{eq: a2.33}\begin{aligned}&\int|u|^{2}(\eta+1)^{2}+\lambda^{2m} \mbox{Re}\;\sum_{|\alpha|=|\beta|=m} w^{-1}a_{\alpha, \beta}\partial^{\beta}(u(\eta+1))\cdot\overline{\partial^{\alpha}(u(\eta+1))}d\upsilon\\
&\quad\quad=\frac{\lambda^{2m}}{2}\int_{\rr^{n+1}}|u|^{2}\partial_{t}(\eta^{2})d\upsilon+\int_{\rr^{n+1}}|u|^{2}(2\eta+1))d\upsilon\\
&\quad\quad+\lambda^{2m} \mbox{Re}\;\int_{\rr^{n+1}}w^{-1}a_{\alpha, \beta}\l[\partial^{\beta}(u(\eta+1))\cdot\overline{\partial^{\alpha}(u(\eta+1))}-\partial^{\beta}u \cdot\overline{\partial^{\alpha}(u(\eta+1)^{2})}\r]d\upsilon\\
&\quad\quad+\lambda^{2m} \mbox{Re}\; \int_{\rr^{n+1}}w^{-1}a_{\alpha, \beta}\partial^{\beta}u \cdot\overline{\partial^{\alpha}(u(2\eta+1))}d\upsilon\\
&:=J_{1}+J_{2}+J_{3}+J_{4}.
\end{aligned}\end{equation} 

We first treat $J_{1}.$ By the definitions of $\tilde{\eta}$ (\eqref{eq: a2.30}) and $\eta,$ also Young's inequality, it is not hard to show 
\begin{equation*}\begin{aligned}
J_{1}& \lesssim \ez \int_{\rr^{n+1}}|u|^{2}\eta^{2}d\upsilon+C(\ez)\lambda^{4m}\int_{\rr^{n+1}}|u|^{2}(\partial_{t}\eta)^{2}d\upsilon\\
&\lesssim \ez \int_{\rr^{n+1}}|u|^{2}\eta^{2}d\upsilon+C(\ez)\lambda^{4m}\Delta^{2}d^{-4m}\int_{\rr^{n+1}}|\eta+1|^{2}|u|^{2}d\upsilon\quad\mbox{$(\lambda<d)$}\\
&\lesssim\ez \int_{\rr^{n+1}}|u|^{2}(\eta+1)^{2}d\upsilon+\kappa^{2}C(\ez)\int_{\rr^{n+1}}|u|^{2}(\eta+1)^{2}d\upsilon+\int_{\rr^{n+1}}|u|^{2}d\upsilon.
\end{aligned}\end{equation*} Next, we turn to the domination for $J_{3}.$ Apparently, using Leibniz's rule, we can write 

\begin{equation*}\begin{aligned}J_{3}&=-\lambda^{2m} \mbox{Re}\;\sum_{|\alpha|=|\beta|=m} \int_{\rr^{n+1}}w^{-1}a_{\alpha, \beta}\sum_{\tau<\alpha}C_{\alpha}^{\tau}\partial^{\beta}u\overline{\partial^{\tau}u} \partial^{\alpha-\tau}(\eta+1)^{2}d\upsilon\\
&\quad\quad+ \lambda^{2m} \mbox{Re}\;\sum_{|\alpha|=|\beta|=m} \int_{\rr^{n+1}}w^{-1}a_{\alpha, \beta}\sum_{|\tau|+|\gamma|<2m}C_{\alpha}^{\tau}C_{\beta}^{\gamma}\overline{\partial^{\tau}u}\partial^{\gamma}u\partial^{\alpha-\tau}(\eta+1)\partial^{\beta-\gamma}(\eta+1)d\upsilon\\
&:=-J_{31}+J_{32}. 
\end{aligned}\end{equation*} 
Moreover, the term $J_{32}$ can be further decomposed into \begin{equation*}\begin{aligned}J_{32}&=\lambda^{2m} \mbox{Re}\;\sum_{|\alpha|=|\beta|=m}\sum_{\tau=\alpha, \gamma<\beta} C_{\beta}^{\gamma}\int_{\rr^{n+1}}w^{-1}a_{\alpha, \beta}\overline{\partial^{\alpha}u}\partial^{\gamma}u(\eta+1)\partial^{\beta-\gamma}(\eta+1)d\upsilon\\
&\quad\quad+\lambda^{2m} \mbox{Re}\;\sum_{|\alpha|=|\beta|=m}\sum_{\tau<\alpha, \gamma=\beta}C_{\alpha}^{\tau} \int_{\rr^{n+1}}w^{-1}a_{\alpha, \beta}\overline{\partial^{\tau}u}\partial^{\beta}u\partial^{\alpha-\tau}(\eta+1)(\eta+1)d\upsilon\\
&\quad\quad+\lambda^{2m} \mbox{Re}\;\sum_{|\alpha|=|\beta|=m} \sum_{\tau<\alpha, \gamma<\beta}\int_{\rr^{n+1}}w^{-1}a_{\alpha, \beta}C_{\alpha}^{\tau}C_{\beta}^{\gamma}\overline{\partial^{\tau}u}\partial^{\gamma}u\partial^{\alpha-\tau}(\eta+1)\partial^{\beta-\gamma}(\eta+1)d\upsilon\\
&:=J_{321}+J_{322}+J_{323}. 
\end{aligned}\end{equation*} 
Observe that $J_{321}$ and $J_{322}$ are essentially the same type by the symmetry of $\alpha, \beta.$ Hence, for simplicity, we only handle $J_{321}$ here. 

Again, by the definition of $\eta$ and Leibniz's rule, for any $|\xi|\leq m,$  \begin{equation}\label{eq: a2.36}\partial^{\xi}(\eta+1)=(\eta+1)P^{\Delta}_{\xi}(\partial_{1},...,\partial_{n})\tilde{\eta}, \end{equation}where $P^{\Delta}_{\xi}$ denotes a homogeneous polynomial of degree $|\xi|$ ($P^{\Delta}_{0}:=1$) satisfying \begin{equation}\label{eq: a2.34}|P^{\Delta}_{\xi}(\partial_{1},...,\partial_{n})\tilde{\eta}|\lesssim \l(\frac{\Delta}{d}\r)^{|\xi|}\quad (\Delta\geq 1).\end{equation} 
By bundling up $u$ and $\eta+1,$ a very tedious calculation leads to  \begin{equation}\label{eq: a2.35}\partial^{\xi}u(\eta+1)=\sum_{\tau\leq \xi}P_{\xi-\tau}^{\Delta}(\partial_{1},...,\partial_{n})\tilde{\eta} \partial^{\tau}(u(\eta+1)).\end{equation} 
Inserting \eqref{eq: a2.36} and \eqref{eq: a2.35} into $J_{321}$ we get
\begin{equation*}\begin{aligned}J_{321}&=\lambda^{2m} \mbox{Re}\;\sum_{|\alpha|=|\beta|=m}\sum_{\tau=\alpha, \gamma<\beta}C_{\beta}^{\gamma}\int_{\rr^{n+1}}w^{-1}a_{\alpha, \beta}\overline{\partial^{\alpha}u}\partial^{\gamma}u(\eta+1)^{2}P^{\Delta}_{\beta-\gamma}(\partial_{1},...,\partial_{n})\tilde{\eta}d\upsilon\\
&=\lambda^{2m} \mbox{Re}\;\sum_{|\alpha|=|\beta|=m}\sum_{\gamma<\beta} \sum_{\tau\leq \alpha}\sum_{S\leq \gamma}C_{\beta}^{\gamma}\int_{\rr^{n+1}}w^{-1}a_{\alpha, \beta}P^{\Delta}_{\beta-\gamma}(\partial_{1},...,\partial_{n})\tilde{\eta}\\ 
&\quad\quad\quad \times  P_{\alpha-\tau}^{\Delta}(\partial_{1},...,\partial_{n})\tilde{\eta}\overline{\partial^{\tau}(u(\eta+1))}P_{\gamma-S}^{\Delta}(\partial_{1},...,\partial_{n})\tilde{\eta}\partial^{S}(u(\eta+1))d\upsilon\\
\end{aligned}\end{equation*}  Since $\gamma<\beta,$ that is, $|S|+|\tau|\leq 2m-1,$ then, by \eqref{eq: a2.34} and \eqref{eq: a2.29},
\begin{equation}\label{eq: a2.37}\begin{aligned}
&\lambda^{2m}\int_{\rr^{n+1}}w^{-1}a_{\alpha, \beta}P^{\Delta}_{\beta-\gamma}(\partial_{1},...,\partial_{n})\tilde{\eta}\\ 
&\quad\quad\times  P_{\alpha-\tau}^{\Delta}(\partial_{1},...,\partial_{n})\tilde{\eta}\overline{\partial^{\tau}(u(\eta+1))}P_{\gamma-S}^{\Delta}(\partial_{1},...,\partial_{n})\tilde{\eta}\partial^{S}(u(\eta+1))d\upsilon\\
&\lesssim\lambda^{2m}c_{2}\l(\frac{\Delta}{d}\r)^{|\alpha-\tau|+|\beta-S|}\|\partial^{S}(u(\eta+1))\|_{2, \upsilon}\|\partial^{\tau}(u(\eta+1))\|_{2, \upsilon}\\
\end{aligned}\end{equation} 
\begin{equation*}\begin{aligned}
&\lesssim c_{2} \kappa \l(\lambda^{|\tau|}\|\partial^{\tau}(u(\eta+1))\|_{2, \upsilon}\r)\l(\lambda^{|S|}\|\partial^{S}(u(\eta+1))\|_{2, \upsilon}\r)\quad (\kappa<1)\\
&\lesssim \kappa c_{2}c_{0}(n, m, [w]_{A_{2}})  \lambda^{|\tau|}\l(\int_{\rr^{n+1}}|u(\eta+1)|^{2}w\r)^{ \frac{(1-\frac{|\tau|}{m})}{2} }\l(\int_{\rr^{n+1}}|\nabla^{m}(u(\eta+1)|^{2})w\r)^{\frac{|\tau|}{2m}}\\
&\quad\quad\times \lambda^{|S|}\l(\int_{\rr^{n+1}}|u(\eta+1)|^{2}w\r)^{\frac{ (1-\frac{|S|}{m})}{2} }\l(\int_{\rr^{n+1}}|\nabla^{m}(u(\eta+1)|^{2})w\r)^{\frac{|S|}{2m}}\\
&\lesssim  \kappa c_{2} c_{0}(n, m, [w]_{A_{2}})\l(\int_{\rr^{n+1}}|u(\eta+1)|^{2}d\upsilon+\lambda^{2m}\int_{\rr^{n+1}}|\nabla^{m}(u(\eta+1))|^{2}d\upsilon\r)
\end{aligned}\end{equation*} where in the second last step we used \eqref{eq: a2.32} and in the last step used Young's inequality!
Remark that $c_{0}$ is independent of $\Delta.$ Armed with \eqref{eq: a2.37}, it follows that $$J_{321}\leq C(m, n, c_{2}, [w]_{A_{2}})\kappa\l(\int_{\rr^{n+1}}|u(\eta+1)|^{2}d\upsilon+\lambda^{2m}\int_{\rr^{n+1}}|\nabla^{m}(u(\eta+1))|^{2}d\upsilon\r).$$

 Similarly, 
 \begin{equation}\label{eq: a2.38}\begin{aligned}
\lambda^{2m}&\int_{\rr^{n+1}}w^{-1}a_{\alpha, \beta}\overline{\partial^{\tau}u}\partial^{\gamma}u\partial^{\alpha-\tau}(\eta+1)\partial^{\beta-\gamma}(\eta+1)d\upsilon\\
&\quad\quad=\lambda^{2m}\int_{\rr^{n+1}}w^{-1}a_{\alpha, \beta}\overline{\partial^{\tau}u}\partial^{\gamma}u(\eta+1)^{2}P^{\Delta}_{\alpha-\tau}(\partial_{1},...,\partial_{n})\tilde{\eta}P^{\Delta}_{\beta-\gamma}(\partial_{1},...,\partial_{n})\tilde{\eta}d\upsilon\\
\end{aligned}\end{equation} 
\begin{equation*}\begin{aligned}
&=\lambda^{2m}\sum_{S\leq \tau}\sum_{\xi\leq \gamma}\int_{\rr^{n+1}}w^{-1}a_{\alpha, \beta}P_{\gamma-\xi}^{\Delta}(\partial_{1},...,\partial_{n})\tilde{\eta} \partial^{\xi}(u(\eta+1))\\
&\quad\quad\times P_{\tau-S}^{\Delta}(\partial_{1},...,\partial_{n})\tilde{\eta} \overline{\partial^{S}(u(\eta+1))}P^{\Delta}_{\beta-\gamma}(\partial_{1},...,\partial_{n})\tilde{\eta}P^{\Delta}_{\alpha-\tau}(\partial_{1},...,\partial_{n})\tilde{\eta}d\upsilon\\
&\lesssim\lambda^{2m}c_{2}\sum_{S\leq \tau}\sum_{\xi\leq \gamma}\l(\frac{\Delta}{d}\r)^{|\alpha-S|+|\beta-\xi|}\|\partial^{\xi}(u(\eta+1))\|_{2, \upsilon}\|\partial^{S}(u(\eta+1))\|_{2, \upsilon}\\
&\lesssim c_{2} \kappa^{2} \sum_{S\leq \tau}\sum_{\xi\leq \gamma}\l(\lambda^{|S|}\|\partial^{S}(u(\eta+1))\|_{2, \upsilon}\r)\l(\lambda^{|\xi|}\|\partial^{\xi}(u(\eta+1))\|_{2, \upsilon}\r)\quad (|\xi|+|S|\leq 2m-2)\\
&\lesssim c_{2} \kappa^{2}\sum_{S\leq \tau}\sum_{\xi\leq \gamma}C(\xi, S, m) \lambda^{|S|}\l(\int_{\rr^{n+1}}|u(\eta+1)|^{2}w\r)^{\frac{ (1-\frac{|S|}{m}) } {2}}\l(\int_{\rr^{n+1}}|\nabla^{m}(u(\eta+1)|^{2})w\r)^{\frac{|S|}{2m}}\\
&\quad\quad\times \lambda^{|\xi|}\l(\int_{\rr^{n+1}}|u(\eta+1)|^{2}w\r)^{\frac{ (1-\frac{|\xi|}{m}) } {2}}\l(\int_{\rr^{n+1}}|\nabla^{m}(u(\eta+1)|^{2})w\r)^{\frac{|\xi|}{2m}}\\
&\lesssim c_{2} \kappa^{2}\l(\int_{\rr^{n+1}}|u(\eta+1)|^{2}d\upsilon+\lambda^{2m}\int_{\rr^{n+1}}|\nabla^{m}(u(\eta+1))|^{2}d\upsilon\r).
\end{aligned}\end{equation*}  Thus, $$J_{323}\leq C(m, n, c_{2}, [w]_{A_{2}})\kappa^{2}\l(\int_{\rr^{n+1}}|u(\eta+1)|^{2}d\upsilon+\lambda^{2m}\int_{\rr^{n+1}}|\nabla^{m}(u(\eta+1))|^{2}d\upsilon\r).$$ Applying a similar argument in \eqref{eq: a2.37}-\eqref{eq: a2.38} to $J_{31}$ we can prove that $J_{31}$ has the same bound as $J_{321}$. The details are omitted. 

It remains to bound $J_{2}$ and $J_{4}.$ At this time, we just exploit the trivial bound $$\|\eta\|_{\infty}\lesssim e^{\Delta}$$ and Young's inequality and to deduce $$J_{2}\lesssim e^{\Delta} \int_{\rr^{n+1}}|u|^{2}d\upsilon $$ and $$J_{4}\lesssim  \lambda^{2m}\int_{\rr^{n+1}}|\nabla^{m}u|^{2}d\upsilon+\ez \lambda^{2m}\int_{\rr^{n+1}}|\nabla^{m}(u(\eta+1))|^{2}d\upsilon.$$ Collecting all the above estimates and letting $\ez$ and $\kappa$ small enough, it follows from \eqref{eq: a2.33} and \eqref{eq: a2.28} that \begin{equation*}\begin{aligned}\int_{\rr^{n+1}}|u|^{2}(\eta+1)^{2}d\upsilon +&\lambda^{2m}\int_{\rr^{n+1}}\sum_{|\alpha|=m}|\partial^{\alpha}(u(\eta+1))|^{2}d\upsilon \\
&\lesssim e^{\Delta} \l(\int_{\rr^{n+1}}|u|^{2}d\upsilon+\lambda^{2m}\int_{\rr^{n+1}}|\nabla^{m}u|^{2}d\upsilon\r).\end{aligned}\end{equation*} Furthermore, invoking Lemma \ref{lemma: a12.7}, we arrive at \begin{equation}\label{eq: a2.39}\int_{\rr^{n+1}}|u|^{2}(\eta+1)^{2}d\upsilon+\lambda^{2m}\int_{\rr^{n+1}}\sum_{|\alpha|=m}|\partial^{\alpha}(u(\eta+1))|^{2}d\upsilon\lesssim e^{\Delta} \int_{\rr^{n+1}}|f|^{2}d\upsilon.\end{equation}
As $\tilde{\eta}=1$ on $F$ and $\eta+1=e^{\Delta \tilde{\eta}},$ then we instantly conclude by \eqref{eq: a2.39} that $$e^{2\Delta}\int_{F}|u|^{2}d\upsilon+e^{2\Delta}\lambda^{2m}\int_{F}\sum_{|\alpha|=m}|\partial^{\alpha}u|^{2}d\upsilon\lesssim e^{\Delta} \int_{E}|f|^{2}d\upsilon.$$ Consequently, \eqref{eq: a2.26} is proved.

Finally, we come to the proof of \eqref{eq: a2.261}. In fact, by a duality argument, \eqref{eq: a2.261} can be attributed to the estimate \eqref{eq: a2.26} for $\Ez_{\lambda}^{*}$. To be precise, by interchanging the roles of $E$ and $F,$  we can derive 

\begin{equation*}\begin{aligned}\int_{F}|\lambda^{m} \Ez_{\lambda}[w^{-1}\partial^{\alpha}(wf)]|^{2}d\upsilon&=\sup_{g}\l(\int_{\rdm}\lambda^{m} \Ez_{\lambda}[w^{-1}\partial^{\alpha}(wf)]\overline{g}d\upsilon\r)^{2}\\
 &=\sup_{g}\l(\int_{E}\lambda^{m} f \overline{\partial^{\alpha} \Ez_{\lambda}^{*} g}d\upsilon\r)^{2}\\
 &\lesssim e^{-\frac{ d}{c\lambda}}\int_{E}|f|^{2}d\upsilon,\end{aligned}\end{equation*} where the supremum is taken with respect to all $g\in L^{2}_{\upsilon}$ with $\supp g\subset F.$
This suffices.

\section{Proof of \eqref{eq: a2.54}: Part II }
As customary in the field, we need to decompose
\begin{equation*}\begin{aligned}
\lambda^{m}\H \Ez_{\lambda} f&=\lambda^{m}\Ez_{\lambda}\H(I-\P_{\lambda}+\P_{\lambda})^{m} f\\
&=\sum_{k=1}^{m}C_{m}^{k}\lambda^{m}\Ez_{\lambda}\H\l(\P_{\lambda}^{k}(I-\P_{\lambda})^{m-k}\r)f+\lambda^{m}\Ez_{\lambda}\H(I-\P_{\lambda})^{m}f
&:=Y_{1}+Y_{2}\end{aligned}\end{equation*} in order to proceed with the proof of \eqref{eq: a2.54}, where $\P_{\lambda}$ denotes the identity approximation operator. This compels us to introduce the parabolic weighted Littlewood-Paley theory of higher order.

\noindent{\textbf{5.1. Parabolic weighted Littlewood-Paley theory of higher order.}} Given two locally integrable functions $h_{1}$ and $h_{2},$ defined on $\rz$ and $\rr$ respectively, the maximal operators associated with them are defined by $$\M^{(1)}(h_{1})(x):=\sup_{r>0}\fint_{Q_{r}(x)}|h_{1}(y)|dy$$ and 
 $$\M^{(2)}(h_{2})(x):=\sup_{r>0}\fint_{I_{r}(t)}|h_{2}(s)|ds.$$ In the sequel, we let $\P_{\tiny{\textcircled{1}}}(x)\in C_{0}^{\infty}(\rz)$ and $\P_{\tiny{\textcircled{2}}}(t)\in C_{0}^{\infty}(\rr)$ be two radial functions, both of which have integral 1. For all $t\in \rr,$ $x\in \rz$ and $\lambda>0,$ we set $$(\P_{\tiny{\textcircled{1}}})_{\lambda}(x):=\lambda^{-n}\P_{\tiny{\textcircled{1}}}(x/\lambda),\quad (\P_{\tiny{\textcircled{2}}})_{\lambda}(t):=\lambda^{-2m}\P_{\tiny{\textcircled{2}}}(t/\lambda^{2m}),$$ and $$\P(t, x):=\P_{\tiny{\textcircled{1}}}(x)\P_{\tiny{\textcircled{2}}}(t),\quad \P_{\lambda}(t, x):=(\P_{\tiny{\textcircled{1}}})_{\lambda}(x)(\P_{\tiny{\textcircled{2}}})_{\lambda}(t).$$ With a slight abuse of notation, we use $\P_{\lambda}$ to represent the convolution operator $$\P_{\lambda}f(t, x):=\P_{\lambda}*f(t, x)=\int_{\rdm}\P_{\lambda}(t-s, x-y)f(y, s)dyds.$$ The same rule applies to $(\P_{\tiny{\textcircled{1}}})_{\lambda}$ and $(\P_{\tiny{\textcircled{2}}})_{\lambda}.$ It is well-known that $$(\P_{\tiny{\textcircled{1}}})_{\lambda}f(x, t)\leq \M^{(1)}(f(\cdot, t))(x),$$
 $$(\P_{\tiny{\textcircled{2}}})_{\lambda}f(x, t)\leq \M^{(2)}(f(x, \cdot))(x)$$ and $$\P_{\lambda}f(x, t)\leq \M^{(1)}(\M^{(2)}f(\cdot, t))(x)$$ almost everywhere for any $f\in L^{2}(\rdm, d\upsilon)$ ($\subset L^{1}_{loc}(\rdm)$ thanks to \eqref{eq: h5}). Moreover, both $\M^{(1)}$ and $\M^{(2)}$ are bounded on $L^{2}(\rdm, d\upsilon)$ by keeping one of the variables fiexed, see \cite{EMS}. Consequently, \begin{equation}\label{eq: a2.52}\sup_{\lambda>0}\l(\|\P_{\lambda}f\|_{2, \upsilon}+\|(\P_{\tiny{\textcircled{2}}})_{\lambda}f\|_{2, \upsilon}+\|(\P_{\tiny{\textcircled{1}}})_{\lambda}f\|_{2, \upsilon}\r)\lesssim \|f\|_{2, \upsilon}.\end{equation}

An argument similar to the one in \cite[Lemma 5.1]{AMK} leads to the following lemma.  

\begin{lemma}\label{lemma: a2.6}\;  For all $f\in L^{2}_{\upsilon}(\rdm),$ then for any $1\leq |\alpha|\leq m,$\begin{equation}\label{eq: a2.40}
|||\lambda^{|\alpha|}\nabla^{|\alpha|}\P_{\lambda}f|||_{2, \upsilon}+|||\lambda^{m}D_{t}^{1/2}\P_{\lambda}f|||_{2, \upsilon}+|||\lambda^{2m}\partial_{t}\P_{\lambda}f|||_{2, \upsilon}\lesssim \|f\|_{2, \upsilon}.\end{equation}  \end{lemma} 
{\it Proof.}\quad Using Plancherel's theorem in time because of the separation between variables $x$ and $t$, 
and the uniform $L^{2}(\rdm, d\upsilon)-$boundedness of $(\P_{\tiny{\textcircled{1}}})_{\lambda}$ and $(\P_{\tiny{\textcircled{2}}})_{\lambda}$ in \eqref{eq: a2.52}, we have \begin{equation*}\begin{aligned}
|||\lambda^{m}D_{t}^{1/2}\P_{\lambda}f|||_{2, \upsilon}^{2}&=\int_{0}^{\infty}\int_{\rdm}|(\P_{\tiny{\textcircled{1}}})_{\lambda}\lambda^{m}D_{t}^{1/2}(\P_{\tiny{\textcircled{2}}})_{\lambda}f|^{2}\frac{d\upsilon d\lambda}{\lambda}\\
&\lesssim \int_{0}^{\infty}\int_{\rdm}|\lambda^{m} D_{t}^{1/2}(\P_{\tiny{\textcircled{2}}})_{\lambda}f|^{2}\frac{d\tau d\lambda}{\lambda}dw(x)\\
&\approx \int_{0}^{\infty}\int_{\rdm}|\lambda^{m}|\tau|^{1/2}\widehat{\P_{\tiny{\textcircled{2}}}}(\lambda^{2m}\tau)\hat{f}(x, \tau)|^{2}\frac{d\tau d\lambda}{\lambda}dw(x)\\
&\lesssim \int_{\rdm}|\hat{f}(x, \tau)|^{2}d\tau dw(x)\int_{0}^{\infty}|\lambda^{m}|\tau|^{1/2}\widehat{\P_{\tiny{\textcircled{2}}}}(\lambda^{2m}\tau)|^{2}\frac{ d\lambda}{\lambda}\\
&\lesssim \int_{\rdm}|f(x, t)|^{2}dt dw(x),\end{aligned}
\end{equation*} where in the last step we used the fact that $\widehat{\P_{\tiny{\textcircled{2}}}}$ is a radial Schwartz function. The same argument applies to $|||\lambda^{2m}\partial_{t}\P_{\lambda}f|||_{2, \upsilon}.$ 

Interchanging the roles of $(\P_{\tiny{\textcircled{1}}})_{\lambda}$ and $(\P_{\tiny{\textcircled{2}}})_{\lambda}$ we also have 

\begin{equation*}\begin{aligned}
|||\lambda^{|\alpha|}\nabla^{|\alpha|}\P_{\lambda}f|||_{2, \upsilon}^{2}&=\int_{0}^{\infty}\int_{\rdm}|(\P_{\tiny{\textcircled{2}}})_{\lambda}\lambda^{|\alpha|}\nabla^{|\alpha|}(\P_{\tiny{\textcircled{1}}})_{\lambda}f|^{2}\frac{d\upsilon d\lambda}{\lambda}\\
&\lesssim \int_{0}^{\infty}\int_{\rdm}|\lambda^{|\alpha|}\nabla^{|\alpha|}(\P_{\tiny{\textcircled{1}}})_{\lambda} f|^{2}\frac{d\tau d\lambda}{\lambda}dw(x)\\
&\lesssim \int_{0}^{\infty}\int_{\rdm}|\lambda^{|\alpha|}\nabla^{|\alpha|}(\P_{\tiny{\textcircled{1}}})_{\lambda}f|^{2}\frac{d\tau d\lambda}{\lambda}dw(x)\\
&\lesssim \int_{\rdm}|f(x, \tau)|^{2}d\tau dw(x),\end{aligned}
\end{equation*} where in the last step we used Lemma \ref{lemma: b1.1} since $\nabla^{|\alpha|}\P_{\tiny{\textcircled{1}}}$ ($|\alpha|\geq 1$) is a Schwartz function such that $\widehat{\nabla^{|\alpha|}\P_{\tiny{\textcircled{1}}}}(0)=0.$

\hfill$\Box$

\begin{lemma}\label{lemma: a2.7}\;  For all $f\in \E_{\upsilon}(\rdm),$  \begin{equation}\label{eq: a2.41}
||| \lambda^{-m}(I-\P_{\lambda})^{m}f|||_{2, \upsilon}\lesssim\|\nabla^{m} f\|_{2, \upsilon} + \|D_{t}^{1/2}f\|_{2, \upsilon}\approx \|\de f\|_{2, \upsilon}.\end{equation}  \end{lemma} 
{\it Proof.}\quad Since $\widehat{\P_{\tiny{\textcircled{2}}}}(0)=1,$ it is easy to see $$|1-\widehat{\P_{\tiny{\textcircled{2}}}}(\lambda^{2m}\tau)|\lesssim \min\{1, \lambda^{2m}|\tau|\}.$$ Thus 

\begin{equation*}\begin{aligned}
\int_{0}^{\infty}|1-\widehat{\P_{\tiny{\textcircled{2}}}}(\lambda^{2m}\tau)|^{2}\frac{d\lambda}{\lambda^{2m+1}}&\lesssim\int_{0}^{\infty}| \min\{1, \lambda^{2m}|\tau|\}|^{2}\frac{d\lambda}{\lambda^{2m+1}}\\
&\lesssim \int_{0}^{|\tau|^{-\frac{1}{2m}}}\lambda^{2m-1}|\tau|^{2}d\lambda+\int_{|\tau|^{-\frac{1}{2m}}}^{\infty}\lambda^{2m+1}d\lambda\\
&\lesssim |\tau|,\end{aligned} 
\end{equation*} which induces $$||| \lambda^{-m}(I-(\P_{\tiny{\textcircled{2}}})_{\lambda})f|||_{2, \upsilon}\lesssim \|D_{t}^{1/2}f\|_{2, \upsilon}.$$
Note that $(\P_{\tiny{\textcircled{1}}})_{\lambda}$ commutes with $(\P_{\tiny{\textcircled{2}}})_{\lambda}.$ Then we have\begin{equation}\label{eq: a2.65}I-\P_{\lambda}=(\P_{\tiny{\textcircled{2}}})_{\lambda}(I-(\P_{\tiny{\textcircled{1}}})_{\lambda})+(I-(\P_{\tiny{\textcircled{2}}})_{\lambda}).\end{equation} 
In addition, it follows from $$(I-\P_{\lambda})^{m}=\sum_{k=0}^{m}C_{m}^{k}[(\P_{\tiny{\textcircled{2}}})_{\lambda}(I-(\P_{\tiny{\textcircled{1}}})_{\lambda})]^{k}(I-(\P_{\tiny{\textcircled{2}}})_{\lambda})^{m-k}$$ that
\begin{equation*}\begin{aligned}
||| \lambda^{-m}(I-\P_{\lambda})^{m}f|||_{2, \upsilon}&\lesssim \sum_{k=0}^{m-1}C_{m}^{k}||| \lambda^{-m}[(\P_{\tiny{\textcircled{2}}})_{\lambda}(I-(\P_{\tiny{\textcircled{1}}})_{\lambda})]^{k}(I-(\P_{\tiny{\textcircled{2}}})_{\lambda})^{m-k}f|||_{2, \upsilon}\\
&\quad\quad +||| \lambda^{-m}[(\P_{\tiny{\textcircled{2}}})_{\lambda}(I-(\P_{\tiny{\textcircled{1}}})_{\lambda})]^{m}f|||_{2, \upsilon}:=G_{1}+G_{2}.
\end{aligned} 
\end{equation*}
Utilizing \eqref{eq: a2.52} again, we get for $m-k\geq 1$ that 
$$G_{1} \lesssim \sum_{k=0}^{m-1}C_{m}^{k}c(k, m)||| \lambda^{-m}(I-(\P_{\tiny{\textcircled{2}}})_{\lambda})f|||_{2, \upsilon}\lesssim \|D_{t}^{1/2}f\|_{2, \upsilon}.$$

Next, we seek to prove, for any $1\leq k\leq m,$ \begin{equation}\label{eq: a2.42}
||| \lambda^{-k}(I-(\P_{\tiny{\textcircled{1}}})_{\lambda})^{k}f|||_{2, \upsilon}\lesssim \|\nabla^{k} f\|_{2, \upsilon}.
\end{equation} Once \eqref{eq: a2.42} is proved, it follows readily that $$G_{2}\lesssim ||| \lambda^{-m}(I-(\P_{\tiny{\textcircled{1}}})_{\lambda})^{m}f|||_{2, \upsilon} \lesssim \|\nabla^{m} f\|_{2, \upsilon}.$$ To show \eqref{eq: a2.42}, we employ the argument below (4.3) in \cite{CR1}. For the sake of presentation, we introduce $$ \widehat{\R_{\lambda}f}:=\l(\frac{(1-\widehat{\P_{\tiny{\textcircled{1}}}}(\lambda \xi)\cdot i\lambda \xi}{|\lambda \xi|^{2}}\r)^{k}\hat{f}:=\widehat{\textbf{R}_{\lambda}*f}.$$ Pick a smooth and radial function $\psi: \rz\to \rr,$ supported in $Q_{1}(0),$ such that \begin{equation}\label{eq: a2.43}\int_{\rz}\xi^{\gamma}\psi(\xi)d\xi=0\quad \mbox{for any}\; |\gamma|\leq m,\end{equation} and $$\int_{0}^{\infty}|\hat{\psi}(\lambda)|^{2}\frac{d\lambda}{\lambda}=1.$$ The existence of such a good function is proven in \cite[Lemma 1.1]{FJW}. We define the associated convolution operator by $$ \Q_{\tau}f:=\psi_{\lambda}*f.$$

Clearly, \begin{equation}\label{eq: a2.48}\lambda^{-k}(I-(\P_{\tiny{\textcircled{1}}})_{\lambda})^{k}f=\R_{\lambda}\nabla^{k}f.\end{equation} Since $\P_{\tiny{\textcircled{1}}}$ is radial and $\nabla \widehat{\P_{\tiny{\textcircled{1}}}}(0)=0,$ it can be derived from \eqref{eq: a2.43} that $$\bigg|\l(\frac{(1-\widehat{\P_{\tiny{\textcircled{1}}}}(\lambda \xi)\cdot i\lambda \xi}{|\lambda \xi|^{2}}\r)^{k}\hat{\psi}(\tau \xi)\bigg| \lesssim \min\{\frac{\lambda^{k}}{\tau^{k}}, \frac{\tau^{k}}{\lambda^{k}}\}\quad \mbox{for any} \; \xi\neq 0, \lambda>0 \;\mbox{and}\; \tau>0.$$
This contributes to \begin{equation}\label{eq: a2.45}\|\R_{\lambda}\Q_{\tau}\|_{L^{2}(\rz)\to L^{2}(\rz)}\lesssim \min\{\frac{\lambda^{k}}{\tau^{k}}, \frac{\tau^{k}}{\lambda^{k}}\}.\end{equation} 
We continue our argument by setting $$\widehat{\textbf{R}_{1}}(\xi):=\frac{(1-\widehat{\P_{\tiny{\textcircled{1}}}}( \xi)\cdot i \xi}{| \xi|^{2}}.$$ According to the equation (35) in \cite[Chapter 4]{AT}, $\textbf{R}_{1}(x)$ admits the following bound: $$|\textbf{R}_{1}(x)|\lesssim \frac{1}{|x|^{n-1}(1+|x|)^{2}}\in L^{1}(\rz).$$ Then, by the fact that $$\textbf{R}(x)=\textbf{R}_{1}(*\textbf{R}_{1})^{k-1}(x),$$ we conclude $$|\textbf{R}(x)|\lesssim \frac{1}{|x|^{n-1}(1+|x|)^{2}}\in L^{1}(\rz).$$
An application of Lemma \ref{lemma: b1.2} yields \begin{equation}\label{eq: a2.44}\sup_{\lambda>0}\|\R_{\lambda}f\|_{L^{2}(w)\to L^{2}(w)}\lesssim 1.\end{equation} By \eqref{eq: a2.44} and the definition of $\Q_{\tau}$, it is straightforward to see that \begin{equation}\label{eq: a2.46}\sup_{\lambda>0, \tau>0}\|\R_{\lambda}\Q_{\tau}\|_{L^{2}(w)\to L^{2}(w)}\lesssim 1.\end{equation} 

Combining \eqref{eq: a2.45} and \eqref{eq: a2.46}, also applying Lemma \ref{lemma: b1.3}, one can obtain that there exists a $0<\theta<1$ such that 
\begin{equation}\label{eq: a2.47}\sup_{\lambda>0, \tau>0}\|\R_{\lambda}\Q_{\tau}\|_{L^{2}(w)\to L^{2}(w)}\lesssim \min\{\frac{\lambda^{k}}{\tau^{k}}, \frac{\tau^{k}}{\lambda^{k}}\}^{\theta}.\end{equation} 
In the rest of the proof, we use Lemma \ref{lemma: b1.4} and \eqref{eq: a2.47} to prove \eqref{eq: a2.41}. The proof can be seen as a variant of that in \cite[Proposition 4.7]{CR}. Fix $t\in \rr$ and $f\in L^{2}_{\upsilon}(\rdm),$ and let $f_{j}$ be defined as in Lemma \ref{lemma: b1.4}. From \eqref{eq: a2.44}, it follows that $$\int_{\rz}|\R_{\lambda}f(t, x)|^{2}w(x)=\lim_{j\to \infty}\int_{\rz}|\R_{\lambda}f_{j}(t, x)|^{2}w(x).$$ Clearly, $$|\R_{\lambda}f_{j}(t, x)|\lesssim \int_{1/j}^{j}|\R_{\lambda}\Q_{\tau}(\chi_{B_{j}}\Q_{\tau}f)(t, x)|\frac{d\tau}{\tau}$$ due to the fact that $\R_{\lambda}$ is sublinear. Exploiting successively Fatou's lemma, Minkowski's inequality and the latter inequality, we arrive at 
\begin{equation*}\begin{aligned}
\int_{0}^{\infty}\int_{\rdm}&|\R_{\lambda}f(t, x)|^{2}w(x)\frac{dtdxd\lambda}{\lambda}\leq \liminf_{j\to \infty}\int_{0}^{\infty}\int_{\rdm}|\R_{\lambda}f_{j}(t, x)|^{2}w(x)\frac{dtdxd\lambda}{\lambda}\\
&\lesssim \liminf_{j\to \infty}\int_{\rr}\int_{0}^{\infty}\int_{\rz}|\R_{\lambda}f_{j}(t, x)|^{2}w(x)\frac{dxd\lambda}{\lambda}dt\\
&\lesssim \liminf_{j\to \infty}\int_{\rr}\int_{0}^{\infty}\int_{\rz}\l(\int_{1/j}^{j}|\R_{\lambda}\Q_{\tau}(\chi_{B_{j}}\Q_{\tau}f)(t, x)|\frac{d\tau}{\tau}\r)^{2}w(x)\frac{dxd\lambda}{\lambda}dt\\
&\lesssim \liminf_{j\to \infty}\int_{\rr}\int_{0}^{\infty}\l(\int_{1/j}^{j}\l(\int_{\rz}|\R_{\lambda}\Q_{\tau}(\chi_{B_{j}}\Q_{\tau}f)(t, x)|^{2}w(x)\r)^{1/2}\frac{d\tau}{\tau}\r)^{2}\frac{d\lambda}{\lambda}dt\\
&\lesssim \liminf_{j\to \infty}\int_{\rr}\int_{0}^{\infty}\l(\int_{1/j}^{j}\min\{\frac{\lambda^{k}}{\tau^{k}}, \frac{\tau^{k}}{\lambda^{k}}\}^{\theta}\|\chi_{B_{j}}\Q_{\tau}f(t, \cdot)\|_{L^{2}(w)}\frac{d\tau}{\tau}\r)^{2}\frac{d\lambda}{\lambda}dt\\
&\lesssim \int_{\rr}\int_{0}^{\infty}\l(\int_{0}^{\infty}\min\{\frac{\lambda^{k}}{\tau^{k}}, \frac{\tau^{k}}{\lambda^{k}}\}^{\theta}\|\Q_{\tau}f(t, \cdot)\|_{L^{2}(w)}\frac{d\tau}{\tau}\r)^{2}\frac{d\lambda}{\lambda}dt:=G.
\end{aligned}
\end{equation*} Note that $$\sup_{\lambda>0}\int_{0}^{\infty}\min\{\frac{\lambda^{k}}{\tau^{k}}, \frac{\tau^{k}}{\lambda^{k}}\}^{\theta}\frac{d\tau}{\tau}\leq C(k, \theta),$$ and the same estimate also holds if we reverse the roles of $\lambda$ and $\tau.$ By this inequality, Schwartz's inequality, Fubini's theorem and Lemma \ref{lemma: b1.1}, we obtain
\begin{equation*}\begin{aligned}
G&\lesssim \int_{\rr}\int_{0}^{\infty} \l(\int_{0}^{\infty}\min\{\frac{\lambda^{k}}{\tau^{k}}, \frac{\tau^{k}}{\lambda^{k}}\}^{\theta}\frac{d\tau}{\tau}\r)\l(\int_{0}^{\infty}\min\{\frac{\lambda^{k}}{\tau^{k}}, \frac{\tau^{k}}{\lambda^{k}}\}^{\theta}\|\Q_{\tau}f(t, \cdot)\|_{L^{2}(w)}^{2}\frac{d\tau}{\tau}\r)\frac{d\lambda}{\lambda}dt\\
&\lesssim  \int_{\rr}\int_{0}^{\infty} \int_{0}^{\infty}\min\{\frac{\lambda^{k}}{\tau^{k}}, \frac{\tau^{k}}{\lambda^{k}}\}^{\theta}\|\Q_{\tau}f(t, \cdot)\|_{L^{2}(w)}^{2}\frac{d\tau}{\tau}\frac{d\lambda}{\lambda}dt\\ 
&\lesssim \int_{\rr}\int_{0}^{\infty}\|\Q_{\tau}f(t, \cdot)\|_{L^{2}(w)}^{2}\int_{0}^{\infty}\min\{\frac{\lambda^{k}}{\tau^{k}}, \frac{\tau^{k}}{\lambda^{k}}\}^{\theta}\frac{d\lambda}{\lambda}dt\frac{d\tau}{\tau}\\
&\lesssim \int_{\rr}\int_{0}^{\infty}\|\Q_{\tau}f(t, \cdot)\|_{L^{2}(w)}^{2}dt\frac{d\tau}{\tau}\\
&\lesssim \int_{\rr}\|f(t, \cdot)\|_{L^{2}(w)}^{2}.\end{aligned} 
\end{equation*} 
Going back to \eqref{eq: a2.48}, we eventually reach \eqref{eq: a2.42}, thereby completing the proof. 

\hfill$\Box$ 

It follows easily from the definition of $\Ez_{\lambda}$ that $$\lambda^{m}\H\Ez_{\lambda}=\lambda^{-m}(I-\Ez_{\lambda}).$$ By this, and employing Lemma \ref{lemma: a12.7} and Lemma \ref{lemma: a2.7}, we have $$|||Y_{2} |||_{2, \upsilon}\lesssim |||\lambda^{-m}(I-\P_{\lambda})^{m}f |||_{2, \upsilon}\lesssim\|\de f\|_{2, \upsilon}.$$ It therefore remains to handle the term $Y_{1}.$ To the end, we introduce the averaging operator and the associated estimates. 

Recall that $\Delta:=I\times Q$ denotes a typical parabolic cube in $\rdm.$

\begin{definition}\label{definition: a2.8}\;(\cite[Definition 5.3]{AMK})\;  We let $\A^{(1)}_{\lambda}, \A^{(2)}_{\lambda}, \A_{\lambda}$ be the dyadic averaging operators in $x, t$ and $(t, x)$ with respect to parabolic scaling, that is, if $\Delta:=I\times Q$ is the (unique) dyadic parabolic cube with $l(\Delta)/2<\lambda<l(\Delta)$ containing $(t, x),$ then $$\A^{(1)}_{\lambda}f(t, x):=\fint_{Q}f(t, y)dy,$$  $$\A^{(2)}_{\lambda}f(t, x):=\fint_{I}f(s, x)ds,$$ and 
$$\A_{\lambda}f(t, x):=\fint_{I\times Q}f(s, y)dsdy=\A^{(2)}_{\lambda}\A^{(1)}_{\lambda}f(t, x)=\A^{(1)}_{\lambda}\A^{(2)}_{\lambda}f(t, x).$$  
\end{definition} It is apparent that 
\begin{equation}\label{eq: a2.50}
\sup_{\lambda>0}\l(\|\A_{\lambda}\|_{L^{2}_{\upsilon}\to L^{2}_{\upsilon}}+\|\A^{(1)}_{\lambda}\|_{L^{2}_{\upsilon}\to L^{2}_{\upsilon}}+\|\A^{(2)}_{\lambda}\|_{L^{2}_{\upsilon}\to L^{2}_{\upsilon}}\r)\lesssim 1
\end{equation} holds by the definitions of the averaging operators and the fact that $\M^{(1)}, \M^{(2)}$ are uniformly bounded on $L^{2}_{\upsilon}.$ 

\begin{lemma}\label{lemma: a2.9}\;  For all $f\in L^{2}_{\upsilon}(\rdm),$  \begin{equation}\label{eq: a2.51}
||| (\A_{\lambda}-\P_{\lambda})f|||_{2, \upsilon}+||| (\A_{\lambda}-\P_{\lambda})\A_{\lambda}f|||_{2, \upsilon}\lesssim \| f\|_{2, \upsilon}.\end{equation}  \end{lemma} 
{\it Proof.}\quad We start with the proof of \begin{equation*}
||| (\A_{\lambda}-\P_{\lambda})f|||_{2, \upsilon}\lesssim \| f\|_{2, \upsilon}.\end{equation*} To the end, we first note that $\A^{(2)}_{\lambda}$ commutes with $(\P_{\tiny{\textcircled{1}}})_{\lambda}.$ Then, $\A_{\lambda}-\P_{\lambda}$ can be split into 
 $$\A_{\lambda}-\P_{\lambda}=\A^{(2)}_{\lambda}(\A^{(1)}_{\lambda}-(\P_{\tiny{\textcircled{1}}})_{\lambda})+(\P_{\tiny{\textcircled{1}}})_{\lambda}(\A^{(2)}_{\lambda}-(\P_{\tiny{\textcircled{2}}})_{\lambda}).$$ Using \eqref{eq: a2.52} and \eqref{eq: a2.50}, we can deduce 

\begin{equation*}\begin{aligned}
||| (\A_{\lambda}-\P_{\lambda})f|||_{2, \upsilon}^{2}&\lesssim \int_{\rr}\int_{0}^{\infty} \|(\A^{(1)}_{\lambda}-(\P_{\tiny{\textcircled{1}}})_{\lambda})f(\cdot, t)\|_{L^{2}(w)}^{2}\frac{d\lambda}{\lambda}dt\\
&\quad\quad+ \int_{\rz}\int_{0}^{\infty} \|(\A^{(2)}_{\lambda}-(\P_{\tiny{\textcircled{2}}})_{\lambda})f(\cdot, x)\|_{L^{2}(\rr; dt)}^{2}\frac{d\lambda}{\lambda}w(x)dx\\
&:=S_{1}+S_{2}.
\end{aligned} 
\end{equation*} 
The first term $S_{1}$ has been addressed in \cite[Lemma 5.4]{AMK}, where the proof utilizes \cite[Lemma 5.2]{CR}. By making a change of variables $\tilde{\lambda}=\lambda^{2m}$ in the integrand $$\int_{0}^{\infty} \|(\A^{(2)}_{\lambda}-(\P_{\tiny{\textcircled{2}}})_{\lambda})f(\cdot, x)\|_{L^{2}(\rr; dt)}^{2}\frac{d\lambda}{\lambda},$$ and revisiting the definition $\A^{(2)}_{\lambda},$ the term $S_{2}$ can be annihilated by directly invoking the unweighetd one-dimensional version of \cite[Lemma 5.2]{CR}. 

Observing that $\P_{\lambda}^{2}$ has the same properties as $\P_{\lambda},$ thus we also have $$||| (\A_{\lambda}-\P_{\lambda}^{2})f|||_{2, \upsilon}\lesssim \| f\|_{2, \upsilon}.$$ It is evident that $$\P_{\lambda}\A_{\lambda}-\A_{\lambda}=\P_{\lambda}(\A_{\lambda}-\P_{\lambda})-(\A_{\lambda}-\P_{\lambda}^{2}),$$ then, by \eqref{eq: a2.52} and $\A_{\lambda}^{2}=\A_{\lambda},$ \begin{equation*}\begin{aligned} 
||| (\A_{\lambda}-\P_{\lambda})\A_{\lambda}f|||_{2, \upsilon}&\lesssim |||\P_{\lambda}(\A_{\lambda}-\P_{\lambda})f|||_{2, \upsilon}+|||(\A_{\lambda}-\P_{\lambda}^{2})f|||_{2, \upsilon}\\
&\lesssim  |||(\A_{\lambda}-\P_{\lambda})f|||_{2, \upsilon}+\| f\|_{2, \upsilon}\\
&\lesssim \| f\|_{2, \upsilon}.\end{aligned} 
\end{equation*} Lemma \ref{lemma: a2.9} is proved.

\hfill$\Box$ 

With the above preparations, we are now in a position to commence the proof of $Y_{1}.$ For this purpose, we first rewrite $Y_{1}$ as 
\begin{equation*}\begin{aligned}Y_{1}&=\sum_{k=1}^{m}C_{m}^{k}\lambda^{m}\Ez_{\lambda}\H\l(\P_{\lambda}^{k}(I-\P_{\lambda})^{m-k}\r)f\\
&=\sum_{k=1}^{m}\sum_{l=0}^{m-k}C_{m}^{k}C_{m-k}^{l}(-1)^{l}\lambda^{m}\Ez_{\lambda}\H(\P_{\lambda}^{k+l}f).\\
\end{aligned}
\end{equation*}
Consequently, it suffices to arrange $$Y_{11}:=\lambda^{m}\Ez_{\lambda}\H(\P_{\lambda}^{j}f)\quad \mbox{for any}\; 1\leq j\leq m.$$
Recall \eqref{eq: b005}, that is, $$\H=D_{t}^{1/2}H_{t}D_{t}^{1/2}+\sum_{|\alpha|=|\beta|=m}(-1)^{|\alpha|}w^{-1}(\partial^{\alpha}w)(w^{-1}a_{\alpha, \beta}\partial^{\beta}).$$ Inserting the equation into $Y_{11}$ we achieve $$Y_{11}=(-1)^{m}\sum_{|\alpha|=|\beta|=m}\amalg_{\alpha, \lambda}(w^{-1}a_{\alpha, \beta}\partial^{\beta}\P_{\lambda}^{j}f)+\lambda^{m}\Ez_{\lambda}D_{t}^{1/2}H_{t}D_{t}^{1/2}\P_{\lambda}^{j}f,$$ where $$\amalg_{\alpha, \lambda}:=(-1)^{m}\lambda^{m}\Ez_{\lambda}w^{-1}\partial^{\alpha}(w\cdot), \quad \forall\; |\alpha|=m.$$ 
Similar to the estimate for $Y_{2}$, by Lemma \ref{lemma: a12.7} and Lemma \ref{lemma: a2.6}, we can deduce 
\begin{equation*}\begin{aligned}
 |||\lambda^{m}\Ez_{\lambda}D_{t}^{1/2}H_{t}D_{t}^{1/2}\P_{\lambda}^{j}f |||_{2, \upsilon}&=|||\lambda^{m}\Ez_{\lambda}\P_{\lambda}^{j-1}D_{t}^{1/2}H_{t}D_{t}^{1/2}\P_{\lambda}f |||_{2, \upsilon}\quad \mbox{$(j\geq 1)$}\\
&\lesssim |||\lambda^{m}D_{t}^{1/2}\P_{\lambda}H_{t}D_{t}^{1/2}f |||_{2, \upsilon}\\
&\lesssim \|D_{t}^{1/2}f\|_{2, \upsilon}.\end{aligned}
\end{equation*} Thus we are left to manage $$Y_{111}:=\amalg_{\alpha, \lambda}(w^{-1}a_{\alpha, \beta}\partial^{\beta}\P_{\lambda}^{j}f)\quad \forall j\geq 1, \;|\alpha|=|\beta|=m. $$

\noindent{\textbf{5.2. Principal part approximation.}} Before starting the proof, we show by the off-diagonal estimate for $\amalg_{\alpha, \lambda}$ constructed in Lemma \ref{lemma: a2.5} that the operator $\amalg_{\alpha, \lambda}$ can be defined on $L^{\infty}(\rdm),$ and introduce the principal part approximation operator, together with the related estimates.

For any parabolic cube $\Delta:=I\times Q\subset \rdm$ with $Q:=Q_{r}(x)\subset \rz$ and $I:=I_{r}(t)\subset \rr$ defined in Section 2.1, we let $$E_{k}(\Delta):=2^{k+1}\Delta\setminus 2^{k+1}\Delta\quad (k=1,2,...), E_{0}(\Delta):=2\Delta.$$
Obviously, from Lemma \ref{lemma: a12.7}, it follows that, for any $|\alpha|= m,$$$\|\amalg_{\alpha, \lambda}f\|_{2, \upsilon}\lesssim \|f\|_{2, \upsilon}.$$ 

\begin{definition}\label{definition: a2.10}\; For $b\in L^{\infty}(\rdm),$ then, for any $|\alpha|= m,$ \begin{equation}\label{eq: a2.59}\amalg_{\alpha, \lambda}b:=\lim_{k\to \infty}\amalg_{\alpha, \lambda}(b1_{2^{k}\Delta}),\end{equation} is well-defined with the limit taking in $L^{2}_{loc, \upsilon}$ and $\Delta$ being any parabolic cube in $\rdm$.
\end{definition} 

\begin{remark}\label{remark: s0}\; We give a brief justification for the reasonableness of the definition \eqref{eq: a2.59}. Let $\Delta'$ be another parabolic cube. Then, there exist two large integers $\tilde{N}, M$ such that $M>\tilde{N}$ and $\Delta'\subset 2^{\tilde{N}-1}\Delta.$ Consequently, by exploiting Lemma \ref{lemma: a2.5}, we deduce 
\begin{equation*}\begin{aligned}
\|\amalg_{\alpha, \lambda}(b1_{2^{M}\Delta\setminus 2^{\tilde{N}}\Delta})\|_{L^{2}_{\upsilon}(\Delta')} &\lesssim \sum_{j=\tilde{N}}^{M-1}\|\amalg_{\alpha, \lambda}(b1_{E_{j}(\Delta)})\|_{L^{2}_{\upsilon}(\Delta')}\\
&\lesssim \sum_{j=\tilde{N}}^{M-1} e^{-\frac{l(\Delta)2^{j-1}}{\lambda}}\|b\|_{L^{\infty}}\|1_{2^{j+1}\Delta}\|_{L^{2}_{\upsilon}(\Delta')}\\
&\lesssim \upsilon(\Delta)^{1/2}\|b\|_{L^{\infty}}\sum_{j=\tilde{N}}^{M-1} e^{-\frac{l(\Delta)2^{j-1}}{\lambda}} (2^{2m}D)^{j+1}:=G_{M, \tilde{N}},
\end{aligned}
\end{equation*} where we have also used the fact that $d\upsilon$ is a doubling measure with the doubling constant $2^{2m}D.$
Apparently, $G_{M, \tilde{N}}$ tends to zero as $M, \tilde{N}\to \infty.$ This contributes to that $\{\amalg_{\alpha, \lambda}(b1_{2^{M}\Delta})\}$ is a Cauchy sequence in $L^{2}_{loc, \upsilon}$ and the above definition is rational. 

We choose two parabolic cubes $\Delta_{1}, \Delta_{2}$ and an integer $N_{0}$ large enough such that $\Delta_{1}\subset 2^{N_{0}}\Delta_{2}.$ If $M>\tilde{N}+N_{0},$ then there is a parabolic cube $\tilde{\Delta}$ such that $2^{M}\Delta_{2} \setminus 2^{\tilde{N}}\Delta_{1} \subset 2^{M}\tilde{\Delta}\setminus 2^{\tilde{N}}\tilde{\Delta}.$ An application of the same argument as above yields that
$$\|\amalg_{\alpha, \lambda}(b1_{2^{M}\Delta_{2}})-\amalg_{\alpha, \lambda}(b1_{2^{\tilde{N}}\Delta_{1}})\|_{L^{2}_{\upsilon}(\Delta')} \lesssim \upsilon(\tilde{\Delta})^{1/2}\|b\|_{L^{\infty}}\sum_{j=\tilde{N}}^{M-1} e^{-\frac{l(\tilde{\Delta})2^{j-1}}{\lambda}} (2^{2m}D)^{j+1},$$ which implies that Definition \ref{definition: a2.10} is independent of the choice of
$\Delta.$ In particular, by using the uniform $L^{2}_{\upsilon}-$ boundedness of $\amalg_{\alpha, \lambda},$ we have \begin{equation}\label{eq: a2.60}\begin{aligned}
\|\amalg_{\alpha, \lambda}b\|_{L_{\upsilon}^{2}(\Delta)}&\leq \|\amalg_{\alpha, \lambda}(b1_{2\Delta})\|_{L^{2}(\Delta)}+\|\amalg_{\alpha, \lambda}(b1_{\rdm\setminus 2\Delta})\|_{L^{2}_{\upsilon}(\Delta)}\\
&\lesssim \upsilon(\Delta)^{1/2}\|b\|_{L^{\infty}}\l(1+\sum_{j=1}^{\infty} e^{-\frac{l(\Delta)2^{j-1}}{\lambda}} (2^{2m}D)^{j+1}\r).\end{aligned}
\end{equation}

\end{remark} 

The following lemma is a direct consequence of Definition \ref{definition: a2.10} coupled with Remark \ref{remark: s0}.  

\begin{lemma}\label{lemma: a2.11}\; For any $b\in L^{\infty}(\rdm)$ and $f\in L^{2}_{\upsilon}.$ Then, for any $|\alpha|= m,$
$$\|(\amalg_{\alpha, \lambda}b)\A_{\lambda}f\|_{L^{2}_{\upsilon}}\lesssim \|b\|_{L^{\infty}}\|f\|_{L^{2}_{\upsilon}}.$$
\end{lemma} 
{\it Proof.}\quad Given a parabolic cube $\Delta \subset \rdm$ such that $l(\Delta)/2<\lambda<l(\Delta).$ From the definition of $\A_{\lambda},$ it follows that $\A_{\lambda}f$ is constant on $\Delta.$ Then, by \eqref{eq: a2.60}, we see $$\int_{\Delta}|(\amalg_{\alpha, \lambda}b)\A_{\lambda}f|^{2}d\upsilon\leq \int_{\Delta}|(\amalg_{\alpha, \lambda}b)|^{2}d\upsilon\cdot \fint_{\Delta}|\A_{\lambda}f|^{2}d\upsilon\lesssim \|b\|_{L^{\infty}}^{2}\int_{\Delta}|\A_{\lambda}f|^{2}d\upsilon.$$ Decomposing $\rdm$ into a grid of cubes $\{\Delta_{j}\}$ with $l(\Delta_{j})/2<\lambda<l(\Delta_{j}),$ we obtain $$\|(\amalg_{\alpha, \lambda}b)\A_{\lambda}f\|^{2}_{L^{2}_{\upsilon}}\lesssim\sum_{\Delta_{j}}\int_{\Delta_{j}}|(\amalg_{\alpha, \lambda}b)\A_{\lambda}f|^{2}d\upsilon\lesssim  \|b\|_{L^{\infty}}^{2}\|\A_{\lambda}f\|^{2}_{L^{2}_{\upsilon}}\lesssim \|b\|_{L^{\infty}}^{2}\|f\|^{2}_{L^{2}_{\upsilon}},$$ where in the last step we used \eqref{eq: a2.50}.

\hfill$\Box$ 

In the sequel, we define $$\amalg_{\lambda}^{\alpha, \beta}w^{-1}a_{\alpha, \beta}:=\amalg_{\alpha, \lambda}(w^{-1}a_{\alpha, \beta})$$ and $$\R_{\lambda}^{\alpha, \beta}f:=\amalg_{\alpha, \lambda}(w^{-1}a_{\alpha, \beta}f)-(\amalg_{\lambda}^{\alpha, \beta}w^{-1}a_{\alpha, \beta})\A_{\lambda}f.$$ 
The operator $\R_{\lambda}^{\alpha, \beta}$ can be regarded as an approximation to $\amalg_{\lambda}^{\alpha, \beta}w^{-1}a_{\alpha, \beta}$ (called ``principal part approximation'' in \cite{AMK}). For this operator, we have the following bound.
 \begin{proposition}\label{proposition: a2.13}\; Let $f\in C^{\infty}(\rdm)\cap L^{2}_{\upsilon}.$ Then, for any $|\alpha|=|\beta|= m,$\begin{equation}\label{eq: a2.63}
\|\R_{\lambda}^{\alpha, \beta}f\|_{2, \upsilon}\lesssim \|\lambda \nabla f\|_{2, \upsilon}+\|\lambda^{2m} \partial_{t} f\|_{2, \upsilon}.\end{equation}
\end{proposition} 
The proof of Proposition \ref{proposition: a2.13} relies on the simple fact, as shown in Lemma \ref{lemma: a2.12}.
\begin{lemma}\label{lemma: a2.12}\; Let $f\in C^{\infty}(\rdm)\cap L^{2}_{\upsilon}.$ For any parabolic cube $\Delta$ and non-negative integer $j,$
\begin{equation}\label{eq: a2.62}\int_{E_{j}(\Delta)}|f-(f)_{\Delta}|^{2}d\upsilon \lesssim (j+1)\l(\int_{2^{j+1}\Delta}2^{2j}l(\Delta)^{2}|\nabla f|^{2}+2^{4mj}l(\Delta)^{4m}|\partial_{t}f|^{2}d\upsilon \r).\end{equation} 
\end{lemma} 
{\it Proof.}\quad Note that $$f-(f)_{\Delta}=f-(f)_{2^{j+1}\Delta}+(f)_{2^{j+1}\Delta}-(f)_{2^{j}\Delta}+...+(f)_{2\Delta}-(f)_{\Delta}.$$ Let $h_{Q}(t):=\fint_{Q}f(t, y)dy.$ We rewrite $$f-(f)_{2^{k+1}\Delta}=(f-h_{2^{k+1}Q}(t))+(h-(h)_{2^{k+1}I}.$$ 
An application of \cite[Theorem 15.26]{JK} (the weighted Poincar$\acute{e}$ inequality in $x-$variable) shows that $$\int_{2^{k+1}\Delta}|f-h_{2^{k+1}Q}(t))|^{2}d\upsilon \lesssim (2^{k+1}l(\Delta))^{2}\int_{2^{k+1}\Delta}|\nabla f|^{2}d\upsilon $$ and $$\int_{2^{k+1}\Delta}|h-\fint_{2^{k+1}I}h(t)dt|^{2}\lesssim (2^{2m(k+1)}l(\Delta)^{2m})^{2}\int_{2^{k+1}\Delta}|\partial_{t} f|^{2}d\upsilon, $$ where we have also used \eqref{eq: a2.61} in the process. Connecting the above two inequalities we reach  \eqref{eq: a2.62}.  

\hfill$\Box$ 

\noindent{\textbf{Proof of Proposition \ref{proposition: a2.13}.}}  Fix $(t, x)\in \rdm$ and $\lambda>0.$ Let $\Delta$ be the unique dyadic parabolic cube with $l(\Delta)/2<\lambda<l(\Delta)$ that contains $(t, x).$ Then, according to the definitions of the two operators $\R_{\lambda}^{\alpha, \beta}$ and $\A_{\lambda},$ we see 
$$\R_{\lambda}^{\alpha, \beta}f(t, x)=\amalg_{\alpha, \lambda}(w^{-1}a_{\alpha, \beta}(f-(f)_{\Delta})(t, x).$$
Exploiting Lemma \ref{lemma: a2.5}, Lemma \ref{lemma: a2.12} and the fact that $2^{j+1}\Delta$ intersects at most $2^{(j+1)(n+2m)}$ cubes in the dyadic decomposition of $\rdm,$ we can deduce 
\begin{equation*}\begin{aligned}
\|\R_{\lambda}^{\alpha, \beta}f\|^{2}_{L^{2}_{\upsilon}} &=\sum_{\Delta}\int_{\Delta}|\amalg_{\alpha, \lambda}(w^{-1}a_{\alpha, \beta}(f-(f)_{\Delta}))|^{2}d\upsilon \\
&\lesssim \sum_{\Delta}\l(\sum_{j=0}^{\infty}\l(\int_{\Delta}|\amalg_{\alpha, \lambda}(w^{-1}a_{\alpha, \beta}(f-(f)_{\Delta}))1_{E_{j}(\Delta)})|^{2}d\upsilon \r)^{1/2}\r)^{2}\\
&\lesssim\sum_{\Delta}\l(\sum_{j=0}^{\infty} e^{-\frac{2^{j}}{c}}\l(\int_{E_{j}(\Delta)}|(f-(f)_{\Delta})|^{2}d\upsilon \r)^{1/2}\r)\\
\end{aligned}
\end{equation*} 
\begin{equation*}\begin{aligned}
&\lesssim\sum_{\Delta}\sum_{j=0}^{\infty} e^{-\frac{2^{j}}{c}}(j+1)\l(\int_{2^{j+1}\Delta}2^{2j}l(\Delta)^{2}|\nabla f|^{2}+2^{4mj}l(\Delta)^{4m}|\partial_{t}f|^{2}d\upsilon \r)\\
&\lesssim\sum_{\Delta}\sum_{j=0}^{\infty} e^{-\frac{2^{j}}{c}}(j+1)2^{4mj}\l(\int_{2^{j+1}\Delta}\lambda^{2}|\nabla f|^{2}+\lambda^{4m}|\partial_{t}f|^{2}d\upsilon \r)\\
&\lesssim\sum_{j=0}^{\infty} e^{-\frac{2^{j}}{c}}(j+1)2^{4mj}2^{(j+1)(n+2m)}\l(\int_{\rdm}\lambda^{2}|\nabla f|^{2}+\lambda^{4m}|\partial_{t}f|^{2}d\upsilon \r)\\
&\lesssim \l(\int_{\rdm}\lambda^{2}|\nabla f|^{2}+\lambda^{4m}|\partial_{t}f|^{2}d\upsilon \r).\end{aligned}
\end{equation*} This suffices. 

\hfill$\Box$ 

To control $Y_{111},$ we proceed as follows. Using the fact that $\A_{\lambda}^{2}=\A_{\lambda},$ we have 
\begin{equation}\label{eq: e0}\begin{aligned}
Y_{111}&=\amalg_{\alpha, \lambda}(w^{-1}a_{\alpha, \beta}\P_{\lambda}^{j}\partial^{\beta}f)\\
&=\R_{\lambda}^{\alpha, \beta}(\P_{\lambda}^{j}\partial^{\beta}f)+\amalg_{\alpha, \lambda}(w^{-1}a_{\alpha, \beta})\A_{\lambda}\P_{\lambda}^{j-1}(\P_{\lambda}-\A_{\lambda})\partial^{\beta}f\\
&\quad\quad+\amalg_{\alpha, \lambda}(w^{-1}a_{\alpha, \beta})\A_{\lambda}\P_{\lambda}^{j-1}\A_{\lambda}\partial^{\beta}f\\
&=\R_{\lambda}^{\alpha, \beta}(\P_{\lambda}^{j}\partial^{\beta}f)+\amalg_{\alpha, \lambda}(w^{-1}a_{\alpha, \beta})\A_{\lambda}\P_{\lambda}^{j-1}(\P_{\lambda}-\A_{\lambda})\partial^{\beta}f\\
&\quad\quad+\amalg_{\alpha, \lambda}(w^{-1}a_{\alpha, \beta})\A_{\lambda}\P_{\lambda}^{j-2}(\P_{\lambda}-\A_{\lambda})\A_{\lambda}\partial^{\beta}f\\
&\quad\quad+\amalg_{\alpha, \lambda}(w^{-1}a_{\alpha, \beta})\A_{\lambda}\P_{\lambda}^{j-2}\A_{\lambda}\partial^{\beta}f\\
&...\\
&=\R_{\lambda}^{\alpha, \beta}(\P_{\lambda}^{j}\partial^{\beta}f)+\amalg_{\alpha, \lambda}(w^{-1}a_{\alpha, \beta})\A_{\lambda}\P_{\lambda}^{j-1}(\P_{\lambda}-\A_{\lambda})\partial^{\beta}f\\
&\quad\quad +\sum_{s=0}^{j-2}\amalg_{\alpha, \lambda}(w^{-1}a_{\alpha, \beta})\A_{\lambda}\P_{\lambda}^{s}(\P_{\lambda}-\A_{\lambda})\A_{\lambda}\partial^{\beta}f\quad (\mbox{if}\; j\geq 2)\\
&\quad\quad+\amalg_{\alpha, \lambda}(w^{-1}a_{\alpha, \beta})\A_{\lambda}\partial^{\beta}f.
\end{aligned}
\end{equation}  
First, applying in succession Proposition \ref{proposition: a2.13}, \eqref{eq: a2.52} and \eqref{eq: a2.40} in Lemma \ref{lemma: a2.6}, we obtain 

\begin{equation*}\begin{aligned}
 |||\R_{\lambda}^{\alpha, \beta}\P_{\lambda}^{j}\partial^{\beta}f |||_{2, \upsilon}&\lesssim |||\lambda \nabla \P_{\lambda}^{j}\partial^{\beta}f|||_{2, \upsilon}+\|\lambda^{2m} \partial_{t} \P_{\lambda}^{j}\partial^{\beta}f\|_{2, \upsilon}\quad \mbox{$(j\geq 1)$}\\
&\lesssim |||\lambda \nabla \P_{\lambda}\partial^{\beta}f|||_{2, \upsilon}+\|\lambda^{2m} \partial_{t} \P_{\lambda}\partial^{\beta}f\|_{2, \upsilon}\lesssim \|\partial^{\beta}f\|_{2, \upsilon}.
\end{aligned}
\end{equation*} Second, utilizing Lemma \ref{lemma: a2.11}, \eqref{eq: a2.52} and Lemma \ref{lemma: a2.9}, we get
\begin{equation*}\begin{aligned} |||\amalg_{\alpha, \lambda}(w^{-1}a_{\alpha, \beta})\A_{\lambda}\P_{\lambda}^{j-1}(\P_{\lambda}-\A_{\lambda})\partial^{\beta}f |||_{2, \upsilon}&\lesssim \|w^{-1}a_{\alpha, \beta}\|_{L^{\infty}} |||\P_{\lambda}^{j-1}(\P_{\lambda}-\A_{\lambda})\partial^{\beta}f |||_{2, \upsilon}\\
&\lesssim \|w^{-1}a_{\alpha, \beta}\|_{L^{\infty}}|||(\P_{\lambda}-\A_{\lambda})\partial^{\beta}f |||_{2, \upsilon}\\
&\lesssim \|w^{-1}a_{\alpha, \beta}\|_{L^{\infty}}\|\partial^{\beta}f\|_{2, \upsilon},\end{aligned}
\end{equation*} and, for any $1\leq s\leq j-2,$ \begin{equation*}\begin{aligned} |||\amalg_{\alpha, \lambda}(w^{-1}a_{\alpha, \beta})\A_{\lambda}\P_{\lambda}^{s}(\P_{\lambda}-\A_{\lambda})\A_{\lambda}\partial^{\beta}f |||_{2, \upsilon}&\lesssim \|w^{-1}a_{\alpha, \beta}\|_{L^{\infty}} |||\P_{\lambda}^{s}(\P_{\lambda}-\A_{\lambda})\A_{\lambda}\partial^{\beta}f |||_{2, \upsilon}\\
&\lesssim \|w^{-1}a_{\alpha, \beta}\|_{L^{\infty}}|||(\P_{\lambda}-\A_{\lambda})\A_{\lambda}\partial^{\beta}f |||_{2, \upsilon}\\
&\lesssim \|w^{-1}a_{\alpha, \beta}\|_{L^{\infty}}\|\partial^{\beta}f\|_{2, \upsilon}.\end{aligned}
\end{equation*} 

Gathering the above estimates, we finally arrive at \begin{equation}\label{eq: a2.64}
|||Y_{111}-\sum_{|\alpha|=|\beta|=m}\amalg_{\alpha, \lambda}(w^{-1}a_{\alpha, \beta})\A_{\lambda}\partial^{\beta}f|||_{2, \upsilon}\lesssim \|\de f\|_{2, \upsilon}.\end{equation} Set $$\hbar_{\lambda}(t, x):=(\sum_{|\alpha|=m}\amalg_{\alpha, \lambda}(w^{-1}a_{\alpha, \beta})(t, x))_{|\beta|=m}.$$ Then $$\sum_{|\alpha|=|\beta|=m}\amalg_{\alpha, \lambda}(w^{-1}a_{\alpha, \beta})\A_{\lambda}\partial^{\beta}f=\hbar_{\lambda}\cdot \A_{\lambda}\nabla^{m}f.$$ If we can show that \begin{equation}\label{eq: s1}|||\hbar_{\lambda}\cdot \A_{\lambda}\nabla^{m}f |||_{2, \upsilon}\lesssim \|\nabla^{m}f\|_{2, \upsilon},\end{equation} then it follows from \eqref{eq: a2.64} and \eqref{eq: s1} that $$|||Y_{111}|||\lesssim \|\nabla^{m} f\|_{2, \upsilon} + \|D_{t}^{1/2}f\|_{2, \upsilon},$$ thereby completing the proof of \eqref{eq: a2.54} by invoking the definition of $Y_{1}, Y_{2}, Y_{111}.$ In particular, by \eqref{eq: a2.64}, \eqref{eq: s2} is true. 

\section{Proof of \eqref{eq: a2.54}: Part III}
Through the successive reductions established earlier, our task boils down to proving \eqref{eq: s1}. 
\noindent{\textbf{6.1. Reducing \eqref{eq: s1} to a Carleson estimate.}} 
To verify \eqref{eq: s1}, we need the following two lemmas. The proof of the second lemma is a straightforward adaptation of \cite[Lemma 8.2]{AMK}, hence we omit the details. 

\begin{lemma}\label{lemma: a2.15}\; For all dyadic cubes $\Delta\subset \rdm,$ $$\int_{0}^{l(\Delta)}\int_{\Delta}|\hbar_{\lambda}(t, x)|^{2}\frac{d\upsilon d\lambda}{\lambda}\lesssim \upsilon(\Delta).$$
\end{lemma} 

\begin{lemma}\label{lemma: a2.16}(\cite[Lemma 8.2]{AMK})\; Let $\mu$ be a Borel measure on $\rdm\times \rr^{+}$ such that 
$$\|\mu\|_{\C}:=\sup_{\Delta}\frac{\mu(\Delta\times (0, l(\Delta)])}{\upsilon(\Delta)}<\infty,$$ where the supremum is taken over all dyadic parabolic cubes $\Delta\subset \rdm.$ Then there exists a constant $c_{0,}$ depending only on $n$ and $[w]_{A_{2}},$ such that for any $L^{2}_{\upsilon},$ $$\int_{0}^{\infty}\int_{\rdm}|\A_{\lambda}f|^{2}d\mu(t,x, \lambda)\leq c_{0}\|\mu\|_{\C} \int_{\rdm}|f|^{2}d\upsilon.$$
 \end{lemma} 
Admit Lemma \ref{lemma: a2.15} for the moment. Then \eqref{eq: s1} follows readily. Indeed, letting $$d\mu(t,x, \lambda):=|(\sum_{|\alpha|=m}\amalg_{\alpha, \lambda}(w^{-1}a_{\alpha, \beta}))_{|\beta|=m}|^{2}\frac{d\upsilon d\lambda}{\lambda},$$ 
we conclude by Lemma \ref{lemma: a2.15} that $d\mu(t,x, \lambda)$ is a Carleson measure. Combining this and Lemma \ref{lemma: a2.16}, we see
\begin{equation*}\begin{aligned}|||\hbar_{\lambda}(t, x)\cdot \A_{\lambda} \nabla^{m}f |||^{2}_{2, \upsilon}\lesssim\int_{0}^{\infty}\int_{\rdm}|(\A_{\lambda}\partial^{\beta}f(t, x))_{|\beta|=m}|^{2}d\mu(t,x, \lambda)\lesssim \|\nabla^{m}f\|_{2, \upsilon}^{2}.
\end{aligned}
\end{equation*} Thus, the proof further reduces to building Lemma \ref{lemma: a2.15}.

\noindent{\textbf{6.2. Proving the Carleson estimate by a weighted $Tb-$type argument.}} To accomplish this, we begin by generalizing the weighted $Tb-$type argument displayed in \cite[Section 8.1-8.2]{AMK} to the higher-order case. Accordingly, the first step is to construct appropriate (local) $Tb-$type test functions.

Let $\xi=(\xi_{\beta})_{|\beta|=m} \in (\cc)^{p}$ with $\xi_{\beta}\in \cc$ and $\sum_{|\beta|=m}|\xi_{\beta}|^{2}=1,$ 
and $\chi, \Theta$ be two smooth functions defined on $\rz$ and $\rr,$ respectively, whose values are in $[0, 1].$ 
In particular,
\emph{ $\chi\equiv 1$ on $[-1/2, 1/2]^{n}$ with its support contained in $[-1, 1]^{n},$ and $\Theta \equiv 1$ on $[-1/2^{2m}, 1/2^{2m}]$ with its support contained in $(-1, 1).$}

Fix a parabolic cube $\Delta$ with its center $(t_{\Delta}, x_{\Delta}),$ that is, $\Delta=(t_{\Delta}-\frac{l(\Delta)^{2m}}{2^{2m}}, t_{\Delta}+\frac{l(\Delta)^{2m}}{2^{2m}})\times (x_{\Delta}-\frac{l(\Delta)}{2}, x_{\Delta}+\frac{l(\Delta)}{2}).$ Then we define $$\gamma_{\Delta}(t, x):=\chi(\frac{x-x_{\Delta}}{l(\Delta)})\Theta(\frac{t-t_{\Delta}}{l(\Delta)^{2m}})$$ and $$\digamma_{\Delta}^{\xi}(t, x):=\gamma_{\Delta}(t, x)(\phi_{\Delta}(x)\cdot \overline{\xi}), \quad \phi_{\Delta}(x):=\l(\frac{(x-x_{\Delta})^{\beta}}{\beta!}\r)_{|\beta|=m}.$$ It is evident that  $\digamma_{\Delta}^{\xi}\in \E_{\upsilon}$ and \begin{equation}\label{eq: a2.72}
\nabla^{m}(\phi_{\Delta}(x)\cdot \overline{\xi})\cdot \xi\equiv 1.
\end{equation} For any $0<\ez\ll 1,$ we can define a test function by $$f_{\Delta, \ez}^{\xi}:=(I+(\ez l(\Delta))^{2m}\H)^{-1}\digamma_{\Delta}^{\xi}=\Ez_{\ez l(\Delta)}\digamma_{\Delta}^{\xi}\in \E_{\upsilon}$$ thanks to the definition of $\digamma_{\Delta}^{\xi}$ and Lemma \ref{lemma: a12.7}.

\begin{lemma}\label{lemma: a2.17}\; Let $\xi, f_{\Delta, \ez}^{\xi}$ be defined as above, and $\ez$ be a degree of freedom. Then, 
\begin{equation*}\begin{aligned}
&(i)\quad \|f_{\Delta, \ez}^{\xi}-\digamma_{\Delta}^{\xi}\|_{2, \upsilon}^{2}\lesssim (\ez l(\Delta))^{2m}\upsilon(\Delta),\\
&(ii)\quad \|\nabla^{m}(f_{\Delta, \ez}^{\xi}-\digamma_{\Delta}^{\xi})\|_{2, \upsilon}^{2}+\|D_{t}^{1/2}(f_{\Delta, \ez}^{\xi}-\digamma_{\Delta}^{\xi})\|_{2, \upsilon}^{2}\lesssim \upsilon(\Delta),\\
&(iii)\quad \|\nabla^{m}f_{\Delta, \ez}^{\xi}\|_{2, \upsilon}^{2}+\|D_{t}^{1/2}f_{\Delta, \ez}^{\xi}\|_{2, \upsilon}^{2}\lesssim \upsilon(\Delta).
\end{aligned}
\end{equation*} 
\end{lemma}
{\it Proof.}\quad We first write by \eqref{eq: b005} that \begin{equation}\label{eq: a2.71}\begin{aligned}
f_{\Delta, \ez}^{\xi}-\digamma_{\Delta}^{\xi}&=-(\ez l(\Delta))^{2m}\Ez_{\ez l(\Delta)}\H \digamma_{\Delta}^{\xi}\\
&=-(\ez l(\Delta))^{2m}\Ez_{\ez l(\Delta)}D_{t}^{1/2}H_{t}D_{t}^{1/2}\digamma_{\Delta}^{\xi}\\
&\quad\quad-\sum_{|\alpha|=|\beta|=m}(-1)^{|\alpha|}(\ez l(\Delta))^{2m}\Ez_{\ez l(\Delta)}w^{-1}(\partial^{\alpha}w)(w^{-1}a_{\alpha, \beta}\partial^{\beta}\digamma_{\Delta}^{\xi}).\end{aligned}
\end{equation} 
By Lemma \ref{lemma: a12.7} and \eqref{eq: a2.29}, that is, the uniform $L^{2}_{\upsilon}-$boundedness of $$(\ez l(\Delta))^{m}\Ez_{\ez l(\Delta)}w^{-1}(\partial^{\alpha}w)\quad \mbox{and}\quad (\ez l(\Delta))^{m}\Ez_{\ez l(\Delta)}D_{t}^{1/2},$$ one can show $$\|f_{\Delta, \ez}^{\xi}-\digamma_{\Delta}^{\xi}\|_{2, \upsilon}^{2}\lesssim (\ez l(\Delta))^{2m} \|\de \digamma_{\Delta}^{\xi}\|_{2, \upsilon}^{2}.$$ From the definition of $\digamma_{\Delta}^{\xi}$ and Leibniz's rule, it follows that, for any $|\gamma|\leq m,$ $$|\partial^{\gamma}\digamma_{\Delta}^{\xi}|\lesssim l(\Delta)^{m-|\gamma|}, $$ hence $$\|\nabla^{m}\digamma_{\Delta}^{\xi}\|_{2, \upsilon}^{2}\lesssim \upsilon(\Delta).$$
 An application of Plancherel's theorem in the $t-$variable leads to that 
 $$\|D_{t}^{1/2}\digamma_{\Delta}^{\xi}\|_{2, \upsilon}^{2}\lesssim \upsilon(\Delta).$$ This concludes $(i).$ Equipped with the two estimates above, the conclusion $(iii)$ is now attributed to $(ii).$ To show $(ii),$ it is sufficient to invoke the equation \eqref{eq: a2.71} again and utilize
the uniform $L^{2}_{\upsilon}-$boundedness in Lemma \ref{lemma: a12.7} of $$(\ez l(\Delta))^{2m}\nabla^{m}\Ez_{\ez l(\Delta)}w^{-1}(\partial^{\alpha}w), \quad (\ez l(\Delta))^{2m}D_{t}^{1/2}\Ez_{\ez l(\Delta)}D_{t}^{1/2}$$ and $$(\ez l(\Delta))^{2m}D_{t}^{1/2}\Ez_{\ez l(\Delta)}w^{-1}(\partial^{\alpha}w), \quad (\ez l(\Delta))^{2m}\nabla^{m}\Ez_{\ez l(\Delta)}D_{t}^{1/2}.$$ 
\hfill$\Box$ 

\begin{lemma}\label{lemma: a2.18}\; Let $\Delta:=I\times Q$ be a parabolic dyadic cube and $\xi, f_{\Delta, \ez}^{\xi}$ be defined as in Lemma \ref{lemma: a2.17}. There exist $\ez\in (0, 1),$ depending only on the $n, m, c_{1}, c_{2}$ and $[w]_{A_{2}}$, and a finite set $\W$ of unit vectors in $(\cc)^{p}$ with its cardinality depending on $\ez, m$ and $n,$ such that 
\begin{equation}\label{eq: h10}\begin{aligned}\sup_{\Delta}\frac{1}{\upsilon(\Delta)}&\int_{0}^{l(\Delta)}\int_{\Delta}|\hbar_{\lambda}(t, x)|^{2}\frac{d\upsilon d\lambda}{\lambda}\\
&\lesssim \sum_{\xi \in \W}\sup_{\Delta}\frac{1}{\upsilon(\Delta)}\int_{0}^{l(\Delta)}\int_{\Delta}|\hbar_{\lambda}(t, x)\cdot\A_{\lambda}\nabla^{m}f_{\Delta, \ez}^{\xi}|^{2}\frac{d\upsilon d\lambda}{\lambda},\end{aligned}
\end{equation} where the sumpremum is taken over all dyadic parabolic cubes $\Delta\subset \rdm.$
\end{lemma}
{\it Proof.}\quad Given a unit vector $\xi\in (\cc)^{p},$ we introduce the cone $$\C_{\xi}^{\ez}:=\{v\in (\cc)^{p}: |v-(v \cdot\overline{\xi})\xi|<\ez |v \cdot \overline{\xi}|\}.$$ Clearly, $ (\cc)^{p}$ is covered by a finite number of such cones $\{\C_{\xi}^{\ez}\},$ and the number of cones depends on on $\ez, m$ and $n.$ Fix a cone $\C_{\xi}^{\ez},$ and set $$\Gamma_{\lambda, \xi}^{\ez}(t, x):=1_{\C_{\xi}^{\ez}}(\hbar_{\lambda}(t, x))\hbar_{\lambda}(t, x).$$

\noindent{\textbf{Step 1: Estimate for $f_{\Delta, \ez}^{\xi}.$}} At this moment, our primary focus is on handling the integration $$\int_{\Delta}(1-\nabla^{m}f_{\Delta, \ez}^{\xi}\cdot \xi)dxdt.$$ We start by transforming the integral into $$(1-\nabla^{m}f_{\Delta, \ez}^{\xi}\cdot \xi)=\nabla^{m}g_{\Delta, \ez}^{\xi}\cdot \xi+(1-\nabla^{m}\digamma_{\Delta}^{\xi}\cdot \xi),$$ where $g_{\Delta, \ez}^{\xi}:=\digamma_{\Delta}^{\xi}-f_{\Delta, \ez}^{\xi}.$
By \eqref{eq: a2.72}, a simple calculation shows that \begin{equation*}\begin{aligned}1-\nabla^{m}\digamma_{\Delta}^{\xi}\cdot \xi&=-\sum_{|\beta|=m}\sum_{l\leq \beta-1}C_{\beta}^{l}\partial^{\beta-l}\gamma_{\Delta}(t, x)\partial^{l}(\phi_{\Delta}(x)\cdot \overline{\xi})\xi_{\beta}\quad (\gamma_{\Delta}\equiv 1\;\mbox{on}\; \Delta)\\
&\quad\quad+ 1-\gamma_{\Delta}(t, x)\nabla^{m}(\phi_{\Delta}(x)\cdot \overline{\xi})\cdot \xi=0\quad \mbox{when} \; (t, x)\in \Delta.\end{aligned}
\end{equation*}This yields $$\int_{\Delta}(1-\nabla^{m}\digamma_{\Delta}^{\xi}\cdot \xi)dxdt=0.$$ We now turn to the contribution of $$L:=\int_{\Delta}\nabla^{m}g_{\Delta, \ez}^{\xi}\cdot \xi dxdt.$$ To the end, we choose a smooth function $\psi: \rdm \to [0, 1]$ such that 

\emph{ $\psi\equiv 1$ on $\Delta_{\tau}:=(1-\tau^{2m})I\times (1-\tau)Q$ with its support contained in $\Delta,$ and $$\|\partial^{\alpha}\psi\|_{L^{\infty}}\lesssim \frac{1}{(\tau l(\Delta))^{|\alpha|}}, \quad \|\partial_{t}\psi\|_{_{L^{\infty}}}\lesssim \frac{1}{(\tau l(\Delta))^{2m}}, \quad \forall\; |\alpha|\leq m,$$ } where $\tau\in (0, 1)$ is to be determined later. Obviously, $$L=\int_{\Delta}(1-\psi)\nabla^{m}g_{\Delta, \ez}^{\xi}\cdot \xi dxdt +\int_{\Delta}\psi \nabla^{m}g_{\Delta, \ez}^{\xi}\cdot \xi dxdt:=L_{1}+L_{2}.$$ By $(ii)$ in Lemma \ref{lemma: a2.17} and \eqref{eq: a2.611} for $A_{2}-$weight $d\upsilon^{-1}(t, x):=w^{-1}(x)dxdt,$ we obtain 
\begin{equation*}\begin{aligned}
|L_{1}|&\lesssim \l(\int_{\Delta}|1-\psi|^{2}d\upsilon^{-1}\r)^{1/2}\l(\int_{\Delta}|\nabla^{m}g_{\Delta, \ez}^{\xi}|^{2}d\upsilon\r)^{1/2}\\
&\lesssim \upsilon^{-1}(\Delta\setminus \Delta_{\tau})^{1/2}\upsilon(\Delta)^{1/2}\\
&\lesssim|I|(w^{-1}(\tau Q))^{1/2}(w( Q))^{1/2}\\
&\lesssim|I|\tau^{\eta}(w^{-1}(Q))^{1/2}(w( Q))^{1/2}\\
&\lesssim\tau^{\eta} [w]_{A_{2}}|\Delta|.
\end{aligned}
\end{equation*} Since $\supp \psi\subset \Delta,$ by integrating by parts, we derive $$L_{2}=(-1)^{m}\int_{\rdm}g_{\Delta, \ez}^{\xi}\nabla^{m}\psi\cdot \xi.$$ Using $(i)$ in Lemma \ref{lemma: a2.17} and repeating the agument for $L_{1}$ we conclude 
\begin{equation*}\begin{aligned}
|L_{2}|&\lesssim \l(\int_{\rdm}|\nabla^{m}\psi|^{2}d\upsilon^{-1}\r)^{1/2}\l(\int_{\rdm}|g_{\Delta, \ez}^{\xi}|^{2}d\upsilon\r)^{1/2}\\
&\lesssim  \frac{1}{(\tau l(\Delta))^{m}}\upsilon^{-1}(\Delta\setminus \Delta_{\tau})^{1/2}(\ez l(\Delta))^{m}\upsilon(\Delta)\\
&\lesssim\frac{\ez^{m}}{\tau^{m}}\tau^{\eta} [w]_{A_{2}}|\Delta|\lesssim \frac{\ez^{m}}{\tau^{m}}[w]_{A_{2}}|\Delta|.
\end{aligned}
\end{equation*} By choosing $\tau$ such that $\tau^{\eta}=\frac{\ez^{m}}{\tau^{m}},$ and summarizing the estimates for $L_{1}$ and $L_{2},$ we reach a conclusion that  
\begin{equation}\label{eq: a2.75}
\frac{1}{|\Delta|}\bigg|\int_{\Delta}(1-\nabla^{m}f_{\Delta, \ez}^{\xi}\cdot \xi)\bigg|\lesssim \ez^{\frac{m\eta}{m+\eta}}.
\end{equation} Furthermore, the conclusion $(iii)$ in Lemma \ref{lemma: a2.17} along with Cauchy-Schwarz inequality contribute to
\begin{equation}\label{eq: a2.76}
\frac{1}{|\Delta|}\int_{\Delta}|\nabla^{m}f_{\Delta, \ez}^{\xi}|dxdt\lesssim \frac{1}{|\Delta|}\l(\int_{\Delta}|\nabla^{m}f_{\Delta, \ez}^{\xi}|^{2}d\upsilon\r)^{1/2}\upsilon^{-1}(\Delta)^{1/2}\lesssim [w]_{A_{2}}.
\end{equation}

\noindent{\textbf{Step 2: The choice of $\ez.$}} Armed with \eqref{eq: a2.75} and \eqref{eq: a2.76}, and letting $\ez$ sufficiently small, we conclude that 
$$\frac{1}{|\Delta|}\int_{\Delta}\mbox{Re}\;\;(\nabla^{m}f_{\Delta, \ez}^{\xi}\cdot\xi)dxdt\geq \frac{9}{10}$$ and $$\frac{1}{|\Delta|}\int_{\Delta}|\nabla^{m}f_{\Delta, \ez}^{\xi}|dxdt\leq C_{\infty}$$ hold for some large constant $C_{\infty},$ depending only on the $\ez, n, m, c_{1}, c_{2}$ and $[w]_{A_{2}}$. Below, we implement a well-known stopping time argument rooted in \cite{AHLT} to select a collection $\T^{1}_{\xi}=\{\Delta'\}$ of non-overlapping dyadic subcubes of $\Delta,$ which are maximal with respect to the property that either \begin{equation}\label{eq: a2.77}
\frac{1}{|\Delta'|}\int_{\Delta'}\mbox{Re}\;\;(\nabla^{m}f_{\Delta, \ez}^{\xi}\cdot\xi)dxdt\leq \frac{4}{5},
\end{equation} or \begin{equation}\label{eq: a2.78}
\frac{1}{|\Delta'|}\int_{\Delta'}|\nabla^{m}f_{\Delta, \ez}^{\xi}|dxdt\geq (5\ez)^{-1},\end{equation} holds. Indeed, we subdivide dyadically $\Delta$ and stop the first time either \eqref{eq: a2.77} or \eqref{eq: a2.78} holds. Then $\T^{1}_{\xi}=\{\Delta'\}$  is a disjoint set of prabolic subcubes of $\Delta.$ Let $\T^{2}_{\xi}=\{\Delta''\}$ be the collection of dyadic subcubes of $\Delta$ not contained in any $\Delta'\in \T^{1}_{\xi}.$ Clearly, for any $\Delta''\in \T^{2}_{\xi},$ \begin{equation}\label{eq: a2.79}
\frac{1}{|\Delta''|}\int_{\Delta''}\mbox{Re}\;\;(\nabla^{m}f_{\Delta, \ez}^{\xi}\cdot\xi)dxdt\geq \frac{4}{5}\end{equation} and \begin{equation}\label{eq: a2.80}
\frac{1}{|\Delta''|}\int_{\Delta''}|\nabla^{m}f_{\Delta, \ez}^{\xi}|dxdt\leq (5\ez)^{-1}.
\end{equation} For simplicity, we let $B_{1}$ be the union of the cubes in $\T^{1}_{\xi}$ for which \eqref{eq: a2.77} holds, and $B_{2}$ be the union of the cubes in $\T^{1}_{\xi}$ for which \eqref{eq: a2.78} holds. Then $$\bigg|\bigcup_{\Delta'\in \T^{1}_{\xi}}\Delta'\bigg|\leq |B_{1}|+|B_{2}|.$$ The fact that the cubes in $\T^{1}_{\xi}$ do not overlap yields $$|B_{2}|\leq (5\ez)\int_{\Delta}|\nabla^{m}f_{\Delta, \ez}^{\xi}|dxdt \leq C_{\infty}(5\ez)|\Delta|.$$ Set $b_{\Delta, \ez}^{\xi}:=1-\mbox{Re}\;\;(\nabla^{m}f_{\Delta, \ez}^{\xi}\cdot\xi).$ Then $$|B_{1}|\leq 5\sum_{\Delta'\in B_{1}}\int_{\Delta'}b_{\Delta, \ez}^{\xi}=5\int_{\Delta}b_{\Delta, \ez}^{\xi}dxdt-5\int_{\Delta\setminus B_{1}}b_{\Delta, \ez}^{\xi}dxdt:=H_{1}+H_{2}.$$ It is easy to derive that $$|H_{1}|\leq C_{0} \ez^{\frac{m\eta}{m+\eta}}|\Delta|$$
taking into account of  \eqref{eq: a2.75}. As for the contribution of $H_{2},$ we have 
\begin{equation*}\begin{aligned}
|H_{2}|&\leq 5|\Delta\setminus B_{1}|+5\upsilon^{-1}(\Delta\setminus B_{1})^{1/2}\l(\int_{\Delta}|\nabla^{m}f_{\Delta, \ez}^{\xi}|^{2}d\upsilon\r)^{1/2}\\
&\leq 5|\Delta\setminus B_{1}|+5C\upsilon^{-1}(\Delta\setminus B_{1})^{1/2}\upsilon(\Delta)^{1/2}\\
&\leq 5|\Delta\setminus B_{1}|+5C|\Delta\setminus B_{1}|^{\eta}|\Delta|^{1-\eta}\\
&\leq 5(1+C_{1}\delta^{-\frac{1}{\eta}})|\Delta\setminus B_{1}|+C_{1}\delta^{1-\eta}|\Delta|,
\end{aligned}
\end{equation*} where in the process we have used in succession the Cauchy-Schwartz inequality, $(iii)$ in Lemma \ref{lemma: a2.17}, \eqref{eq: a2.611} and Young's inequality. Summarizing the estimates of $H_{1}$ and $H_{2}$ we arrive at $$|B_{1}|\leq \frac{5+C_{0}\ez^{\frac{m\eta}{m+\eta}}+C_{1}\delta^{1-\eta} +C_{1}\delta^{-\frac{1}{\eta}}}{6+C_{1}\delta^{-\frac{1}{\eta}}}|\Delta|.$$ From the latter inequality, by letting $\ez$ small enough, it can be concluded that there exists a constant $\eta'\in (0, 1),$ depending only on the $\ez, n, m, c_{1}, c_{2}$ and $[w]_{A_{2}}$, such that \begin{equation}\label{eq: a2.81}\bigg|\bigcup_{\Delta'\in \T^{1}_{\xi}}\Delta'\bigg|\leq (1-\eta')|\Delta|.\end{equation}
Owing to \eqref{eq: a2.611} again (see also the argument in \cite[Page 196]{EMS}), it follows instantly from \eqref{eq: a2.81} that \begin{equation}\label{eq: a2.82}\upsilon(\bigcup_{\Delta'\in \T^{1}_{\xi}}\Delta')\leq (1-\eta'')\upsilon(\Delta),\end{equation} where $\eta'' \in (0, 1)$ depends only on the $\ez, n, m, c_{1}, c_{2}$ and $[w]_{A_{2}}$.

\noindent{\textbf{Step 3: Plugging the averaging operator.}} For any $\Delta''\in \T^{2}_{\xi},$ we set $$\vartheta:=\frac{1}{|\Delta''|}\int_{\Delta''}\nabla^{m}f_{\Delta, \ez}^{\xi}\cdot dxdt\in (\cc)^{p}.$$ If $(t, x)\in\Delta''$ and $l(\Delta'')/2<\lambda<l(\Delta''),$ then $\vartheta=\A_{\lambda}\nabla^{m}f_{\Delta, \ez}^{\xi}(t, x)$ by the definition of the averaging operator $\A_{\lambda}.$ Assume that $v:=\hbar_{\lambda}(t, x)\in \C_{\xi}^{\ez}.$ By \eqref{eq: a2.79} and \eqref{eq: a2.80}, it is easy to see that  $$|v-(v\cdot \overline{\xi})\xi|<\ez |v\cdot \overline{\xi}|, \quad \mbox{Re}\; (\vartheta\cdot \xi)\geq \frac{4}{5}, \quad|\vartheta|\leq (5\ez)^{-1}.$$ Revisiting the definition of the inner product on $ (\cc)^{p}$ and adapting the proof of \cite[Lemma 5.10]{AHLT}, we can prove that $$|v|\leq 5|v\cdot \vartheta|.$$ The details are left to the interested reader. Thus, \begin{equation}\label{eq: a2.83}|\Gamma_{\lambda, \xi}^{\ez}(t, x)|\leq 5|\hbar_{\lambda}(t, x)\cdot \A_{\lambda}\nabla^{m}f_{\Delta, \ez}^{\xi}(t, x)|.\end{equation} 

By Whitney decomposition, the Carleson region $\Delta\times (0, l(\Delta)]$ can be partitioned as a union of boxes $\Delta'\times (0, l(\Delta')]$ for $\Delta'\in \T^{1}_{\xi}$ and Whitney rectangles $\Delta''\times (l(\Delta'')/2, l(\Delta'')]$ for $\Delta''\in \T^{2}_{\xi}.$ This allows us to write 
\begin{equation*}\begin{aligned}
\frac{1}{\upsilon(\Delta)}\int_{0}^{l(\Delta)}\int_{\Delta}|\Gamma_{\lambda, \xi}^{\ez}(t, x)|^{2}\frac{d\upsilon d\lambda}{\lambda}&\leq \frac{1}{\mu(\Delta)}\sum_{\Delta'\in \T^{1}_{\xi}}\int_{0}^{l(\Delta')}\int_{\Delta'}|\Gamma_{\lambda, \xi}^{\ez}(t, x)|^{2}\frac{d\upsilon d\lambda}{\lambda}\\
&\quad +\frac{1}{\upsilon(\Delta)}\sum_{\Delta''\in \T^{2}_{\xi}}\int_{l(\Delta'')/2}^{l(\Delta'')}\int_{\Delta''}|\Gamma_{\lambda, \xi}^{\ez}(t, x)|^{2}\frac{d\upsilon d\lambda}{\lambda}\\
&:=Z_{1}+Z_{2}.
\end{aligned}
\end{equation*} To handle $Z_{1},$ we temporarily assume that \begin{equation}\label{eq: a2.84}
\Fz^{\ez}_{\xi}:=\sup_{\tilde{\Delta}}\frac{1}{\upsilon(\tilde{\Delta})}\int_{0}^{l(\tilde{\Delta})}\int_{\tilde{\Delta}}|\Gamma_{\lambda, \xi}^{\ez}(t, x)|^{2}\frac{d\upsilon d\lambda}{\lambda}<\infty\end{equation} is true, where the supremum is taken over all dyadic cubes $\tilde{\Delta}\subset \Delta.$ With \eqref{eq: a2.84} in hand, it follows immediately from \eqref{eq: a2.82} that $$Z_{1}\leq \Fz^{\ez}_{\xi} \frac{1}{\upsilon(\Delta)} \sum_{\Delta' \in \T^{1}_{\xi} } \upsilon(\Delta')\leq (1-\eta'')\Fz^{\ez}_{\xi}.$$ Regarding $Z_{2},$ we can derive from  \eqref{eq: a2.83} that $$Z_{2}\leq \frac{100}{\mu(\Delta)}\int_{0}^{l(\Delta)}\int_{\Delta}|\hbar_{\lambda}(t, x)\cdot \A_{\lambda}\nabla^{m}f_{\Delta, \ez}^{\xi}(t, x)|^{2}\frac{d\mu d\lambda}{\lambda}.$$ Observing that the above estimates hold for all dyadic subcubes of $\Delta,$ we reach a conclusion that $$\Fz^{\ez}_{\xi}\leq (1-\eta'')\Fz^{\ez}_{\xi}+\sup_{\tilde{\Delta}}\frac{100}{\upsilon(\tilde{\Delta})}\int_{0}^{l(\tilde{\Delta})}\int_{\tilde{\Delta}}|\hbar_{\lambda}(t, x)\cdot \A_{\lambda}\nabla^{m}f_{\tilde{\Delta}, \ez}^{\xi}(t, x)|^{2}\frac{d\upsilon d\lambda}{\lambda}.$$ By \eqref{eq: a2.84} again, $$\Fz^{\ez}_{\xi} \lesssim \sup_{\tilde{\Delta}}\frac{1}{\mu(\tilde{\Delta})}\int_{0}^{l(\tilde{\Delta})}\int_{\tilde{\Delta}}|\hbar_{\lambda}(t, x)\cdot \A_{\lambda}\nabla^{m}f_{\tilde{\Delta}, \ez}^{\xi}(t, x)|^{2}\frac{d\mu d\lambda}{\lambda},$$ which implies \eqref{eq: h10}.

To remove the a priori assumption \eqref{eq: a2.84}, we replace $\Gamma_{\lambda, \xi}^{\ez}(t, x)$ by $\Gamma_{\lambda, \xi}^{\ez}(t, x)1_{\{\delta<\lambda<\delta^{-1}\}}$ for $\delta>0$ small in the above process because
\begin{equation*}\begin{aligned}\sup_{\tilde{\Delta}}\frac{1}{\upsilon(\tilde{\Delta})}\int_{0}^{l(\tilde{\Delta})}\int_{\tilde{\Delta}}|1_{\{\delta<\lambda<\delta^{-1}\}}\Gamma_{\lambda, \xi}^{\ez}(t, x)|^{2}\frac{d\upsilon d\lambda}{\lambda}&\lesssim \sup_{\tilde{\Delta}}\int_{\delta}^{l(\tilde{\Delta})}\l(1+\sum_{j=1}^{\infty} e^{-\frac{l(\tilde{\Delta})2^{j-1}}{c\lambda}} (2^{2m}D)^{j+1}\r)^{2}\frac{d\lambda}{\lambda}\\
&\lesssim \int_{1}^{\frac{l(\Delta)}{\delta}}\l(1+\sum_{j=1}^{\infty} e^{-\frac{\lambda 2^{j-1}}{c}} (2^{2m}D)^{j+1}\r)^{2}\frac{d\lambda}{\lambda}<\infty,\end{aligned}  \end{equation*} taking into account of \eqref{eq: a2.60} and $l(\Delta)>l(\tilde{\Delta})>\delta.$ Note that the (implict) constants below \eqref{eq: a2.83} are independent of $\delta.$ Therefore, a limiting argument ($\delta \to 0$) enables us to achieve the desired estimate, which completes the proof of Lemma \ref{lemma: a2.18}. 

\hfill$\Box$ 

We are now ready to prove Lemma \ref{lemma: a2.15}.

\noindent{\textbf{Proof of Lemma \ref{lemma: a2.15}.}}  In veiw of Lemma \ref{lemma: a2.18}, it suffices to show that \begin{equation}\label{eq: a2.85}\sup_{\Delta}\frac{1}{\upsilon(\Delta)}\int_{0}^{l(\Delta)}\int_{\Delta}|\hbar_{\lambda}(t, x)\cdot \A_{\lambda}\nabla^{m}f_{\tilde{\Delta}, \ez}^{\xi}(t, x)|^{2}\frac{d\upsilon d\lambda}{\lambda} \lesssim 1\end{equation} holds for any $\xi \in \W.$ In fact, an application of the triangle inequality leads to
\begin{equation*}\begin{aligned}\frac{1}{\upsilon(\Delta)}\int_{0}^{l(\Delta)}\int_{\Delta}|\hbar_{\lambda}(t, x)\cdot \A_{\lambda}&\nabla^{m}f_{\Delta, \ez}^{\xi}(t, x)|^{2}\frac{d\upsilon d\lambda}{\lambda} \\
&\lesssim \frac{1}{\upsilon(\Delta)}|||(\lambda^{m}\H \Ez_{\lambda} -\hbar_{\lambda}\cdot \A_{\lambda}\nabla^{m})f_{\Delta, \ez}^{\xi}|||_{2, \upsilon}^{2}\\
&\quad\quad +\frac{1}{\upsilon(\Delta)} \int_{0}^{l(\Delta)}\int_{\Delta}|\lambda^{m}\H \Ez_{\lambda} f_{\Delta, \ez}^{\xi}(t, x)|^{2}\frac{d\upsilon d\lambda}{\lambda}\\
&:=D_{1}+D_{2}.
\end{aligned}
\end{equation*} 
Combining \eqref{eq: s2} and Lemma \ref{lemma: a2.17} we conclude $$D_{1}\lesssim \frac{1}{\upsilon(\Delta)}(\|\nabla^{m} f\|_{2, \upsilon} + \|D_{t}^{1/2}f\|_{2, \upsilon}) \lesssim 1.$$ By the trivial equality $$\H f_{\Delta, \ez}^{\xi}=\frac{\digamma_{\Delta}^{\xi}-f_{\Delta, \ez}^{\xi}}{(\ez l(\Delta))^{2m}},$$ Lemma \ref{lemma: a12.7} and Lemma \ref{lemma: a2.17}, we finally arrive at  

\begin{equation*}\begin{aligned}D_{2}&\lesssim \frac{1}{\upsilon(\Delta)} \int_{0}^{l(\Delta)}\lambda^{2m}(\ez l(\Delta))^{-4m}\|\digamma_{\Delta}^{\xi}-f_{\Delta, \ez}^{\xi}\|_{2, \upsilon}^{2}\frac{d\upsilon d\lambda}{\lambda}\\
&\lesssim \frac{1}{\upsilon(\Delta)}\|\digamma_{\Delta}^{\xi}-f_{\Delta, \ez}^{\xi}\|_{2, \upsilon}^{2}(\ez l(\Delta))^{-4m}l(\Delta)^{2m} \\
&\lesssim\ez^{-2m}.
\end{aligned}
\end{equation*} The proof of \eqref{eq: a2.85} is therefore complete. 

\hfill$\Box$ 

At the end, we point out that the proof procedure of Theorem \ref{theorem: a0} applies to Theorem \ref{theorem: a00} with some necessary and obvious modifications. Working out the details along the process of proof for Theorem \ref{theorem: a0} is left to interested readers.

\section{Appendix}

Lemma \ref{lemma: b1.1}-\ref{lemma: b1.4} below are the main technical lemmas developed in \cite{CR, CR1} to deal with the Kato problem for weighted second order elliptic operators. As shown in Section 5 (also in \cite[Section 5]{AMK}), they are indispensable in establishing the weighted Littlewood-Paley theory in the parabolic setting.


\begin{lemma}\label{lemma: b1.1}(\cite[Lemma 4.6]{CR})\; Given $w\in A_{2},$ and let $\psi$ be a Schwartz function such that $\hat{\psi}(0)=0.$ Then for all $f\in L^{2}(w),$ \begin{equation}\label{eq: b1.1}
\int_{0}^{\infty}\int_{\rz}|\psi_{\lambda}*f(x)|^{2}\frac{dw(x)d\lambda}{\lambda}\leq C(n, \psi, [w]_{A_{2}} ) \|f\|^{2}_{L^{2}(w)}.\end{equation}  \end{lemma} 

\begin{lemma}\label{lemma: b1.2}(\cite[Lemma 4.1]{CR})\; Suppose $|\varphi(x)|\leq \Psi(x)$ with $\Psi\in L^{1}(\rz)$ a radial function. Then for any $w\in A_{2}$ the operators $\lambda\to \varphi_{\lambda}*f$ are uniformly bounded on $L^{2}(w).$ Moreover, 
\begin{equation}\label{eq: b1.2}
\sup_{\lambda>0}\|\varphi_{\lambda}*f\|_{L^{2}(w)\to L^{2}(w)}\leq C(n, [w]_{A_{2}} ) \|\Psi\|_{L^{1}(\rz)}.\end{equation}  \end{lemma} 

\begin{lemma}\label{lemma: b1.3}(\cite[Lemma 4.10]{CR})\; Suppose that a sublinear operator $T$ is bounded on $L^{2}(w)$ for all $w\in A_{2}$ with $\|T\|_{L^{2}(w)\to L^{2}(w)}$ depending on $[w]_{A_{2}}$ and the dimension $n.$ Then for any $w\in A_{2},$ there exists a $0<\theta<1$ depending on $[w]_{A_{2}}$ such that $$\|T\|_{L^{2}(w)\to L^{2}(w)}\leq C(n, [w]_{A_{2}}) \|T\|^{1-\theta}_{L^{2}(w^{s})\to L^{2}(w^{s})} \|T\|^{\theta}_{L^{2}(\rz)\to L^{2}(\rz)}$$ where $s>1$ is such that $w^{s}\in A_{2}.$
\end{lemma} 

\begin{lemma}\label{lemma: b1.4}(\cite[Lemma 4.8]{CR})\; For all $w\in A_{2}$ and $f\in L^{2}(w),$ $$\int_{0}^{\infty}\Q_{\tau}^{2}f(x)\frac{d\tau}{\tau}=f(x),$$ where this equality is understood as follows: for each $j>1,$ let $B_{j}$ be the ball centered at 0 of radius $j,$ and define the function $$f_{j}(x)=\int_{1/j}^{j}\Q_{\tau}(\chi_{B_{j}}\Q_{\tau}f)(x)\frac{d\tau}{\tau}.$$ Then for each $j,$ $f_{j}\in L^{2}(w)$ and $\{f_{j}\}$ converges to $f$ in $L^{2}(w).$
\end{lemma} 
 
The following lemma is an analogue of \cite[Lemma 3.3]{AMK}, which plays a part in building the off-diagonal estimates for the resolvent operators. 

\begin{lemma}\label{lemma: b1.5}\; The following statements are true:
 \begin{equation*}\begin{aligned} &(i)\quad \mbox{The space} \; C_{0}^{\infty}(\rdm)\;\mbox{is dense in }\;\E_{\upsilon}.\\ 
 &(ii)\quad \mbox{Multiplication by} \; C^{m}(\rdm)\mbox{-functions is bounded on}\;\E_{\upsilon}.\end{aligned}\end{equation*} 
\end{lemma} 
{\it Proof.}\quad The proof is a straightforward adaptation of the one in \cite[Lemma 3.3]{AMK}. Given $f\in \E_{\upsilon}.$ Without loss of generality, we can assume that $f$ is smooth since convolutions with smooth mollifiers, separately in time and space, provide smooth approximations in $\E_{\upsilon}.$ Thus, it remains to show that $f$ can be approximated by functions in $C_{0}^{\infty}(\rdm)$ in order to conclude $(i).$ 

For this purpose, we choose a sequence $\{\psi_{j}\}\subset C_{0}^{\infty}(\rdm)$ such that $$\|\psi_{j}\|_{L^{\infty}}+j^{|\alpha|}\|\partial^{\alpha}\psi_{j}\|_{L^{\infty}}+j^{2m}\|\partial_{t}\psi_{j}\|_{L^{\infty}}\lesssim 1\quad (\forall\; |\alpha|\leq m)$$ and $\psi_{j} \to 1$ pointwise a.e. as $j\to \infty.$ Set $f_{j}:=\psi_{j}f.$ Obviously, by Leibniz's rule and dominated convergence, we can deduce $$\|\partial^{\alpha}(f_{j})-\partial^{\alpha}f\|_{2, \upsilon}\leq \|\partial^{\alpha}f-\psi_{j} \partial^{\alpha}f\|_{2, \upsilon}+\sum_{1\leq |\gamma|\leq |\alpha|}C_{\alpha}^{\gamma}\|\partial^{\gamma}\psi_{j} \partial^{\alpha-\gamma}f\|_{2, \upsilon}\to 0$$ as $j\to \infty$ for any $|\alpha|\leq m.$ On the other hand, repeating the argument in \cite[(3.3)]{AMK}, we can derive 
\begin{equation*}\begin{aligned}\frac{|(f-f_{j})(t, x)-(f-f_{j})(s, x)|^{2}}{|t-s|^{2}}&\lesssim |f(t, x)|^{2}\frac{|\psi_{j} (t, x)-\psi_{j} (s, x)|^{2}}{|t-s|^{2}}\\
&\quad\quad+\frac{|f(t, x)-f(s, x)|^{2}}{|t-s|^{2}}|\psi_{j}(s, x)-1| \\
&\lesssim \min\{\frac{1}{j^{4m}}, \frac{1}{|t-s|^{2}}\}|f(t, x)|^{2}+\frac{|f(t, x)-f(s, x)|^{2}}{|t-s|^{2}}
\end{aligned}\end{equation*} Recalling the formular \eqref{eq: h2} and the fact $f\in \E_{\upsilon},$ also using the dominated convergence theorem, we get $\|D_{t}^{1/2}(f_{j}-f)\|_{2, \upsilon}\to 0$ as $j\to \infty.$ This completes the proof of $(i).$

In a same fashion, we can conclude, for any $\psi\in C^{m}(\rdm),$ that $$\|\psi f\|_{2, \upsilon}\lesssim \|\psi\|_{L^{\infty}} \|f\|_{2, \upsilon},$$ $$\|\partial^{\alpha}(\psi f)\|_{2, \upsilon}\lesssim \sum_{|\gamma|\leq |\alpha|}C_{\alpha}^{\gamma}\|\partial^{\gamma}\psi\|_{L^{\infty}}\|f\|_{\E_{\upsilon}} \quad (\forall\; |\alpha|\leq m),$$ and 
$$\|D_{t}^{1/2}(\psi f)\|_{2, \upsilon}\lesssim (1+\|\psi\|_{L^{\infty}})^{1/2}\|\partial_{t}\psi\|_{L^{\infty}}^{1/2} \|f\|_{2, \upsilon} +(1+\|\psi\|_{L^{\infty}})\|D_{t}^{1/2}f\|_{2, \upsilon}.$$  This implies $(ii).$

\hfill$\Box$

\section*{Availability of data and material}
 Not applicable.
 
 \section*{Competing interests}
 The author declares that they have no competing interests.



\end{document}